\documentclass[12pt,draft,leqno]{article}
\usepackage{amssymb, eucal, latexsym}

\newtheorem{theorem}{Theorem}[section]
\newtheorem{lm}[theorem]{Lemma}

\newtheorem{cor}[theorem]{Corollary}

\newtheorem{defi}[theorem]{Definition}

\newtheorem{nota}[theorem]{Notation}
\newtheorem{notas}[theorem]{Notations}
\newtheorem{rem}[theorem]{Remark}
\newtheorem{rems}[theorem]{Remarks}
\newtheorem{fact}[theorem]{Fact}

\newtheorem{nist}[theorem]{}

\def\p{\varphi}
\def\a{\alpha}

\def\d{\delta}

\def\g{\gamma}
\def\GA{\Gamma}

\def\l{\lambda}
\def\LAM{\Lambda}

\def\s{\sigma}

\def\pl{\varphi_\Lambda}
\def\fs{\hat{f}}
\def\ps{\hat{\varphi}}

\def\lra{\longrightarrow}

\def\lr{\leftrightarrow}

\def\sbe{\subseteq}
\def\spe{\supseteq}
\def\stm{\setminus}
\def\ems{\emptyset}
\def\nes{\neq\emptyset}

\def\unl{\underline}

\def\fa{\forall}

\def\we{\wedge}

\def\ap{^\prime}
\def\inv{^{-1}}
\def\st{\ |\ }

\def\nin{\not\in}

\def\card #1{\vert #1 \vert}

\def\CC{{\cal C}}

\def\TT{{\cal T}}

\def\K{{\bf K}}
\def\KCon{{\bf KCon}}
\def\Bo{{\bf Bool}}

\def\OHCC{{\bf OpCCon}}
\def\OPLCC{{\bf OpPerLCCon}}
\def\OACC{{\bf DOpCCon}}
\def\OPALC{{\bf DOpPerLCCon}}
\def\SHCC{{\bf SkeCCon}}
\def\SLCC{{\bf SkePerLCCon}}
\def\SACC{{\bf DSkeCCon}}
\def\SALC{{\bf DSkePerLCCon}}

\def\OHC{{\bf OpC}}

\def\OLC{{\bf OpLC}}
\def\OPLC{{\bf OpPerLC}}
\def\OAC{{\bf DOpC}}
\def\OAL{{\bf DOpLC}}
\def\OPAL{{\bf DOpPerLC}}

\def\SAC{{\bf DSkeC}}
\def\SAL{{\bf DSkePerLC}}
\def\SKAL{{\bf DSkeLC}}
\def\SHC{{\bf SkeC}}
\def\SLC{{\bf SkePerLC}}
\def\SKLC{{\bf SkeLC}}

\def\ESAL{{\bf ESkePerLC}}
\def\ESKAL{{\bf ESkeLC}}
\def\EOAL{{\bf EOpLC}}
\def\EOPC{{\bf EOpC}}
\def\EOPAL{{\bf EOpPerLC}}
\def\ESHC{{\bf ESkeC}}
\def\ESKLC{{\bf ESkeLC}}
\def\ESLC{{\bf ESkePerLC}}
\def\EOLC{{\bf EOpLC}}
\def\EOPLC{{\bf EOpPerLC}}

\def\2{\mbox{{\bf 2}}}
\def\3{\mbox{{\bf 3}}}

\def\int{\mbox{{\rm int}}}
\def\Fr{\mbox{{\rm Fr}}}
\def\cl{\mbox{{\rm cl}}}
\def\CL{\mbox{{\rm Clust}}}
\def\BClu{\mbox{{\rm BClust}}}
\def\Ult{\mbox{{\rm Ult}}}

\def\doc{\hspace{-1cm}{\em Proof.}~~}
\def\sq{\hspace*{\fill} \hbox{\vrule\vbox{\hrule\phantom{o}\hrule}\vrule}}
\def\sqs{\sq \vspace{2mm}}

\def\BBBB{{\rm I}\!{\rm B}}

\title{{\LARGE\bf
Some Generalizations of Fedorchuk Duality Theorem -- I}\\
\vspace{0.35cm}
{\large\bf Georgi Dimov}\thanks{This paper was supported by the
project MI-1510/2007 $``$Applied Logics and Topological
Structures" of the Bulgarian Ministry of Education and Science.
}\\
\vspace{0.25cm}
 {\footnotesize Dept. of Math. and
Informatics, Sofia University,  5 J. Bourchier Blvd., 1164 Sofia,
Bulgaria}}

\author{}

\date{}

\begin{document}
\maketitle
\begin{abstract}
{\footnotesize
\noindent Generalizing Duality Theorem of V. V. Fedorchuk
\cite{F}, we prove Stone-type duality theorems for the following
four categories: all of them have as objects the locally compact
Hausdorff spaces, and their morphisms are, respectively, the
continuous skeletal maps, the quasi-open perfect maps, the open
maps, the open perfect maps. In parti\-cular, a Stone-type duality
theorem for the category of all compact Hausdorff spaces and all
open maps between them is proved.  We also obtain equivalence
theorems for these four categories.  The versions of these
theorems for the full subcategories of these categories having as
objects all locally compact connected Hausdorff spaces are
formulated as well.}
\end{abstract}

{\footnotesize {\em  MSC:} primary 54D45, 18A40; secondary 06E15,
54C10, 54E05.

}

\footnotetext[1]{{\footnotesize {\em E-mail address:}
gdimov@fmi.uni-sofia.bg}}

\baselineskip = 0.987\normalbaselineskip

\section*{Introduction}

According to the famous Stone Duality Theorem (\cite{ST}), the
category of all zero-dimensi\-o\-nal compact Hausdorff spaces and
all continuous maps between them is dually equivalent to the
category $\Bo$ of all Boolean algebras and all Boolean
homomorphisms between them. In 1962, H. de Vries \cite{dV}
introduced the notion of {\em compingent Boolean algebra}\/ and
proved that the category of all compact Hausdorff spaces and all
continuous maps between them is dually equivalent to the category
of all complete compingent Boolean algebras and appropriate
morphisms between them. Using de Vries' Theorem, V. V. Fedorchuk
\cite{F} showed that the category $\SHC$ of all compact Hausdorff
spaces and all quasi-open maps between them is dually equivalent
to the category $\SAC$ of all complete normal contact algebras and
all complete Boolean homomorphisms between them satisfying one
simple condition (see Theorem \ref{dcompn} below). The {\em normal
contact algebras}\/ (briefly, NCAs) are Boolean algebras with an
additional relation, called {\em contact relation}. The axioms
which this contact relation satisfies are very similar to the
axioms of Efremovi\v{c} proximities (\cite{EF}). The notion of
normal contact algebra was introduced by Fedorchuk \cite{F} under
the name {\em Boolean $\d$-algebra}\/ as an equivalent expression
of the notion of compingent Boolean algebra of de Vries. We call
such algebras $``$normal contact algebras" because they form a
subclass of the class of {\em contact algebras}\/ introduced in
\cite{DV}. In 1997, Roeper \cite{R} defined the notion of {\em
region-based topology}\/  as one of the possible formalizations of
the ideas of De Laguna \cite{dL} and Whitehead \cite{W} for a
region-based theory of space.  Following \cite{VDDB, DV}, the
region-based topologies of Roeper appear here as {\em local
contact algebras}\/ (briefly, LCAs), because the axioms which they
satisfy almost coincide with the axioms of local proximities of
Leader \cite{LE}. In his paper \cite{R}, Roeper proved the
following theorem: there is a bijective correspondence between all
(up to homeomorphism) locally compact Hausdorff spaces and all (up
to isomorphism) complete LCAs. It generalizes the theorem of de
Vries \cite{dV} that there exists a bijective correspondence
between all (up to homeomorphism) compact Hausdorff spaces and all
(up to isomorphism) complete NCAs. Using the results of Fedorchuk
\cite{F} and Roeper \cite{R}, we show here that the bijective
correspondence established by Roeper can be extended to a duality
between the category $\SKLC$ of all locally compact Hausdorff
spaces and all skeletal (in the sense of Mioduszewski and Rudolf
\cite{MR}) continuous maps between them and the category $\SKAL$
of all complete LCAs and all complete Boolean homomorphisms
between them satisfying two simple axioms; this is done in Theorem
\ref{maintheoremnk} which generalizes the Fedorchuk Duality
Theorem cited above. Further, we regard the non-full subcategory
$\OLC$ (resp., $\OHC$) of the category $\SKLC$ (resp., $\SHC$):
its objects are all locally compact (resp., all compact) Hausdorff
spaces and its morphisms are all open  maps. We find the
corresponding subcategory $\OAL$ (resp., $\OAC$) of the category
$\SKAL$ (resp., $\SAC$) which is dually equivalent to the category
$\OLC$ (resp., $\OHC$) (see Theorem \ref{maintheorem} and Theorem
\ref{dcomp}); as far as we know, even the compact case (i.e. the
result about the category $\OHC$) is new. The subcategories $\SAL$
and $\OPAL$ of the category $\SKAL$ which are dually equivalent,
respectively, to the categories $\SLC$ and $\OPLC$ of all locally
compact Hausdorff spaces and all quasi-open perfect maps
(respectively, all open perfect maps) between them are found as
well (see Theorem \ref{maintheoremn} and Theorem
\ref{maintheoremp}). The versions of all mentioned above theorems
for the  full subcategories of these categories having as objects
all locally compact (resp., compact) connected Hausdorff spaces
are formulated (see Theorems \ref{conthen} and \ref{conthe}).


Following the ideas of Fedorchuk's paper \cite{F}, we define five
categories  $\EOPC$,  $\EOAL$,  $\ESKAL$,  $\ESAL$ and $\EOPAL$,
which are equivalent, respectively, to the categories $\OHC$,
$\OLC$, $\SKLC$,  $\SLC$ and $\OPLC$ (see Theorems \ref{coeth},
\ref{lcoeth}, \ref{lceth}, \ref{lcskpeth}, \ref{lcopeth}). The
equivalence between the categories $\SKLC$ and $\ESKAL$ was almost
established in Roeper's paper \cite{R} (see \ref{Ro} below for
more details). The proof of this equivalence is a slight
modification of the proof of the analogous theorem of Fedorchuk
\cite{F} concerning the case of compact Hausdorff spaces.

Some further development of the results presented here is given in
the second part \cite{D2} of this paper. Let us also mention that
in \cite{D} a category dually equivalent to the category of all
locally compact Hausdorff spaces and all perfect maps between them
is defined, generalizing in this way  de Vries Duality Theorem.

We now fix the notations.

If $\CC$ denotes a category, we write $X\in \card\CC$ if $X$ is
 an object of $\CC$, and $f\in \CC(X,Y)$ if $f$ is a morphism of
 $\CC$ with domain $X$ and codomain $Y$.

All lattices are with top (= unit) and bottom (= zero) elements,
denoted respectively by 1 and 0. We do not require the elements
$0$ and $1$ to be distinct.


Let $X$ and $Y$ be sets. If $f:X\lra Y$ is a function  then for
every subset $Z$ of $Y$, we denote by $f_Z$ the restriction of $f$
with domain $f\inv(Z)$ and codomain $Z$, i.e. $f_Z:f\inv(Z)\lra
Z$.


If $(X,\tau)$ is a topological space and $M$ is a subset of $X$,
we denote by $\cl_{(X,\tau)}(M)$ (or simply by $\cl(M)$ or
$\cl_X(M)$) the closure of $M$ in $(X,\tau)$ and by
$\int_{(X,\tau)}(M)$ (or briefly by $\int(M)$ or $\int_X(M)$) the
interior of $M$ in $(X,\tau)$. The Alexandroff compactification of
a locally compact Hausdorff space $X$ is denoted by $\a X$.

The  closed maps and the open maps between topological spaces are
assumed to be continuous but are not assumed to be onto. Recall
that a map is {\em perfect}\/ if it is closed and compact (i.e.
point inverses are compact sets). A continuous map $f:X\lra Y$ is
{\em irreducible}\/ if $f(X)=Y$ and if, for each proper closed
subset $A$ of $X$, $f(A)\neq Y$.

\section{Preliminaries}
%

\begin{defi}\label{conalg}
\rm
An algebraic system $\underline {B}=(B,0,1,\vee,\we, {}^*, C)$ is
called a {\it contact algebra} (abbreviated as CA) (\cite{DV}) if
$(B,0,1,\vee,\we, {}^*)$ is a Boolean algebra (where the operation
$``$complement" is denoted by $``\ {}^*\ $")
  and $C$
is a binary relation on $B$, satisfying the following axioms:

\smallskip

\noindent (C1) If $a\not=0$ then $aCa$;\\
(C2) If $aCb$ then $a\not=0$ and $b\not=0$;\\
(C3) $aCb$ implies $bCa$;\\
(C4) $aC(b\vee c)$ iff $aCb$ or $aCc$.

\smallskip

\noindent Usually, we shall simply write $(B,C)$ for a contact
algebra. The relation $C$  is called a {\em  contact relation}.
When $B$ is a complete Boolean algebra, we will say that $(B,C)$
is a {\em complete contact algebra}\/ (abbreviated as CCA). For
every  two subsets $M$ and $N$  of $B$, we will write $MCN$ when
$mCn$, for every $m\in M$ and every $n\in N$.

We will say that two CA's $(B_1,C_1)$ and $(B_2,C_2)$ are  {\em
CA-isomorphic} iff there exists a Boolean isomorphism $\p:B_1\lra
B_2$ such that, for each $a,b\in B_1$, $aC_1 b$ iff $\p(a)C_2
\p(b)$. Note that in this paper, by a $``$Boolean isomorphism" we
understand an isomorphism in the category $\Bo$.

A CA $(B,C)$  is called {\em connected}\/ if it satisfies the
following axiom:

\smallskip

\noindent (CON) If $a\neq 0,1$ then $aCa^*$.

\smallskip

A contact algebra $(B,C)$ is called a {\it  normal contact
algebra} (abbreviated as NCA) (\cite{dV,F}) if it satisfies the
following axioms (we will write $``-C$" for $``not\ C$"):

\smallskip

\noindent (C5) If $a(-C)b$ then $a(-C)c$ and $b(-C)c^*$ for some $c\in B$;\\
(C6) If $a\not= 1$ then there exists $b\not= 0$ such that
$b(-C)a$.

\smallskip

\noindent A normal CA is called a {\em complete normal contact
algebra} (abbreviated as CNCA) if it is a CCA.

Note that if $0\neq 1$ then the axiom (C2) follows from the axioms
(C6) and (C4).

For any CA $(B,C)$, we define a binary relation  $``\ll_C $"  on
$B$ (called {\em non-tangential inclusion})  by $``\ a \ll_C b
\leftrightarrow a(-C)b^*\ $". Sometimes we will write simply
$``\ll$" instead of $``\ll_C$".
\end{defi}

The relations $C$ and $\ll$ are inter-definable. For example,
normal contact algebras could be equivalently defined (and exactly
in this way they were defined (under the name of compingent
Boolean algebras) by de Vries in \cite{dV}) as a pair of a Boolean
algebra $B=(B,0,1,\vee,\we,{}^*)$ and a binary relation $\ll$
subject to the following axioms:

\smallskip

\noindent ($\ll$1) $a\ll b$ implies $a\leq b$;\\
($\ll$2) $0\ll 0$;\\
($\ll$3) $a\leq b\ll c\leq t$ implies $a\ll t$;\\
($\ll$4) $a\ll c$ and $b\ll c$ implies $a\vee b\ll c$;\\
($\ll$5) If  $a\ll c$ then $a\ll b\ll c$  for some $b\in B$;\\
($\ll$6) If $a\neq 0$ then there exists $b\neq 0$ such that $b\ll
a$;\\
($\ll$7) $a\ll b$ implies $b^*\ll a^*$.

\smallskip

Note that if $0\neq 1$ then the axiom ($\ll$2) follows from the
axioms ($\ll$3), ($\ll$4), ($\ll$6) and ($\ll$7).

\smallskip

Obviously, contact algebras could be equivalently defined as a
pair of a Boolean algebra $B$ and a binary relation $\ll$ subject
to the  axioms ($\ll$1)-($\ll$4) and ($\ll$7).

\smallskip

It is easy to see that axiom (C5) (resp., (C6)) can be stated
equivalently in the form of ($\ll$5) (resp., ($\ll$6)).

The next notion is a lattice-theoretical counterpart of the
corresponding notion from the theory of proximity spaces (see
\cite{NW}):

\begin{nist}\label{defcluclan}
\rm Let $(B,C)$ be a CA. Then  a non-empty subset $\s $ of $B$ is
called a {\em cluster in} $(B,C)$ (see \cite{VDDB}) if the
following conditions are satisfied:

\smallskip

\noindent (K1) If $a,b\in\s $ then $aCb$;\\
(K2) If $a\vee b\in\s $ then $a\in\s $ or $b\in\s $;\\
(K3) If $aCb$ for every $b\in\s $, then $a\in\s $.

\smallskip

\noindent The set of all clusters in $(B,C)$ will be denoted by
$\CL(B,C)$.

\smallskip

\noindent The set of all ultrafilters in a Boolean algebra $B$
will be denoted by $\Ult(B)$.
\end{nist}

The next three assertions can be proved exactly as Lemma 5.6,
Theorem 5.8 and Corollary 5.10 of \cite{NW}:

\begin{fact}\label{fact29}{\rm (\cite{VDDB})}
If $\s_1,\s_2$ are two clusters in a normal contact algebra
$(B,C)$ and $\s_1\sbe \s_2$ then $\s_1=\s_2$.
\end{fact}

\begin{theorem}\label{conclustth}{\rm (\cite{VDDB})}
A subset $\s$ of a normal contact algebra $(B,C)$ is a cluster iff
there exists an ultrafilter $u$ in $B$ such that
\begin{equation}\label{ultclu}
\s=\{a\in B\st aCb \mbox{ for every } b\in u\}.
\end{equation}
Moreover, given $\s$ and $a_0\in \s$, there exists an ultrafilter
$u$ in $B$ satisfying (\ref{ultclu}) which contains $a_0$.

Note that everywhere in this assertion we can substitute the word
$``$ultrafilter" for $``$basis of an ultrafilter".
\end{theorem}

\begin{cor}\label{uniqult}{\rm (\cite{VDDB})}
Let $(B,C)$ be a normal contact algebra and $u$ be an ultrafilter
(or a basis of an ultrafilter) in $B$. Then there exists a unique
cluster $\s_u$ in $(B,C)$ containing $u$, and
\begin{equation}\label{sigmau}
\s_u=\{a\in B\st aCb \mbox{  for every } b\in u\}.
\end{equation}
\end{cor}

\begin{defi}\label{endb}
\rm In analogy to the corresponding definitions in the theory of
proximity spaces
(see, e.g., \cite{NW}), we say that:\\
(a) a subset $\xi$ of an NCA $(B,C)$ is called an {\em end}\/ if
the following conditions are satisfied:

\smallskip

\noindent (E1)  for any $b,c\in\xi$ there exists  $a\in\xi$ such
that $a\neq 0$,
$a\ll b$ and $a\ll c$;\\
(E2) if $a,b\in B$ and $a\ll b$ then either $a^*\in\xi$ or
$b\in\xi$;

\smallskip

\noindent (b) a subset $v$ of an NCA $(B,C)$ is called a {\em
round filter}\/ if it is a filter and for every $b\in v$ there
exists $a\in v$ such that $a\ll b$.
\end{defi}

The next two theorems (and their proofs) are analogous to the
Theorems 6.7 and 6.11 in \cite{NW} (and their proofs),
respectively:

\begin{theorem}\label{roend}
Let $(B,C)$ be a normal contact algebra and $\xi$ be an end in
$(B,C)$. Then $\xi$ is a maximal round filter in $(B,C)$.
 \end{theorem}

\begin{theorem}\label{cluend}
Let $(B,C)$ be a normal contact algebra and $\s\sbe B$. Then
$\s\in\CL(B,C)$ iff $d(\s)=\{b\in B\st b^*\not\in\s\}$ is an end
in $(B,C)$.
 \end{theorem}

\begin{cor}\label{cluendcor}
Let $(B,C)$ be a normal contact algebra, $\s\in\CL(B,C)$, $a\in B$
and $a\not\in\s$. Then there exists  $b\in B$ such that
$b\not\in\s$ and $a\ll b$.
 \end{cor}

\doc  Put $\xi=d(\s)(=\{c\in B\st c^*\not\in\s\})$.
Then, by \ref{cluend} and \ref{roend}, $\xi$ is a round filter in
$(B,C)$. Since $a\not\in\s$, we obtain that $a^*\in\xi$. Hence,
there exists  $b^*\in\xi$ such that $b^*\ll a^*$. Then
$b\not\in\s$ and $a\ll b$. \sqs

\begin{nist}\label{rct}
\rm Recall that a subset $F$ of a topological space $(X,\tau)$ is
called {\em regular closed}\/ if $F=\cl(\int (F))$. Clearly, $F$
is regular closed iff it is a closure of an open set.

For any topological space $(X,\tau)$, the collection $RC(X,\tau)$
(we will often write simply $RC(X)$) of all regular closed subsets
of $(X,\tau)$ becomes a complete Boolean algebra
$(RC(X,\tau),0,1,\we,\vee,{}^*)$ under the following operations:
$$ 1 = X,  0 = \emptyset, F^* = \cl(X\stm F), F\vee G=F\cup G,
F\we G =\cl(\int(F\cap G)).
$$
The infinite operations are given by the formulas
$\bigvee\{F_\g\st \g\in\GA\}=\cl(\bigcup\{F_\g\st
\g\in\GA\})(=\cl(\bigcup\{\int(F_\g)\st \g\in\GA\})),$ and
$\bigwedge\{F_\g\st \g\in\GA\}=\cl(\int(\bigcap\{F_\g\st
\g\in\GA\})).$

It is easy to see that setting $F \rho_{(X,\tau)} G$ iff $F\cap
G\not = \ems$, we define a contact relation on $RC(X,\tau)$; it is
called a {\em standard contact relation}. So,
$(RC(X,\tau),\rho_{(X,\tau)})$ is a CCA (it is called a {\em
standard contact algebra}). We will often write simply $\rho_X$
instead of $\rho_{(X,\tau)}$. Note that, for $F,G\in RC(X)$,
$F\ll_{\rho_X}G$ iff $F\sbe\int_X(G)$.

Clearly, if $(X,\tau)$ is a normal Hausdorff space then the
standard contact algebra $(RC(X,\tau),\rho_{(X,\tau)})$ is a
complete NCA.

For every topological space $(X,\tau)$, we denote by $RO(X,\tau)$
(or simply by $RO(X)$) the set of all regular open subsets of $X$
(recall that a subset is {\em regular open}\/ if its complement is
regular closed).
\end{nist}

\begin{fact}\label{confact}{\rm (\cite{BG})}
Let $(X,\tau)$ be a topological space. Then the standard contact
algebra $(RC(X,\tau),\rho_{(X,\tau)})$ is connected iff the space
$(X,\tau)$ is connected.
\end{fact}

\begin{nota}\label{sigmax}
\rm
 Let $(X,\tau)$ be a topological space and $x\in X$. Then we
set:
\begin{equation}\label{sxvx}
\s_x=\{F\in RC(X)\st x\in F\} \mbox{ and } \nu_x=\{F\in RC(X)\st
x\in\int(F)\}.
\end{equation}
(Since in our notations the points of a topological space are
denoted only by the letters $``$x,y,z", there will be no confusion
with the notation $\s_u$ introduced in \ref{uniqult}.)
\end{nota}

\begin{fact}\label{sxcluster}
For any topological space $(X,\tau)$ and every point $x\in X$,
$\nu_x$ is a filter in $RC(X)$. If $X$ is regular then $\s_x$ is a
cluster in the CA $(RC(X),\rho_X)$.
\end{fact}

The next notion is a lattice-theoretical counterpart of the
Leader's notion of {\em local proximity} (\cite{LE}):

\begin{defi}\label{locono}{\rm (\cite{R})}
\rm An algebraic system $\underline {B}_{\, l}=(B,0,1,\vee,\we,
{}^*, \rho, \BBBB)$ is called a {\it local contact algebra}
(abbreviated as LCA)   if $(B,0,1, \vee,\we, {}^*)$ is a Boolean
algebra, $\rho$ is a binary relation on $B$ such that $(B,\rho)$
is a CA, and $\BBBB$ is an ideal (possibly non proper) of $B$,
satisfying the following axioms:

\smallskip

\noindent(BC1) If $a\in\BBBB$, $c\in B$ and $a\ll_\rho c$ then
$a\ll_\rho b\ll_\rho c$ for some $b\in\BBBB$  (see \ref{conalg}
for
$``\ll_\rho$");\\
(BC2) If $a\rho b$ then there exists an element $c$ of $\BBBB$
such that
$a\rho (c\we b)$;\\
(BC3) If $a\neq 0$ then there exists  $b\in\BBBB\stm\{0\}$ such
that $b\ll_\rho a$.

\smallskip

Usually, we shall simply  write  $(B, \rho,\BBBB)$ for a local
contact algebra.  We will say that the elements of $\BBBB$ are
{\em bounded} and the elements of $B\stm \BBBB$  are  {\em
unbounded}. When $B$ is a complete Boolean algebra,  the LCA
$(B,\rho,\BBBB)$ is called a {\em complete local contact algebra}
(abbreviated by CLCA).

We will say that two local contact algebras $(B,\rho,\BBBB)$ and
$(B_1,\rho_1,\BBBB_1)$ are  {\em LCA-isomorphic} iff there exists
a Boolean isomorphism $\p:B\lra B_1$ such that, for $a,b\in B$,
$a\rho b$ iff $\p(a)\rho_1 \p(b)$, and $\p(a)\in\BBBB_1$ iff
$a\in\BBBB$.

\end{defi}

\begin{rem}\label{conaln}
\rm Note that if $(B,\rho,\BBBB)$ is a local contact algebra and
$1\in\BBBB$ then $(B,\rho)$ is a normal contact algebra.
Conversely, any normal contact algebra $(B,C)$ can be regarded as
a local contact algebra of the form $(B,C,B)$.
\end{rem}

The following lemmas from \cite{VDDB} are lattice-theoretical
counterparts of some theorems from Leader's paper \cite{LE}.

\begin{lm}\label{Alexprn}{\rm (\cite{VDDB})}
Let $(B,\rho,\BBBB)$ be a local contact algebra. Define a binary
relation $``C_\rho$" on $B$ by
\begin{equation}\label{crho}
aC_\rho b\ \mbox{ iff }\ a\rho b\ \mbox{ or }\ a,b\not\in\BBBB.
\end{equation}
Then $``C_\rho$", called the\/ {\em Alexandroff extension of}\/
$\rho$, is a normal contact relation on $B$ and $(B,C_\rho)$ is a
normal contact algebra.
\end{lm}

\begin{lm}\label{neogrn}{\rm (\cite{VDDB})}
Let $\underline {B}_{\, l}=(B,\rho,\BBBB)$ be a local contact
algebra and let $1\not\in\BBBB$. Then $\s_\infty^{\underline
{B}_{\, l}}=\{b\in B\st b\not\in\BBBB\}$ is a cluster in
$(B,C_\rho)$ (see \ref{Alexprn} for the notation $``C_\rho$").
(Sometimes we will simply  write  $\s_\infty$ or $\s^B_\infty$
instead of $\ \s_\infty^{\underline {B}_{\, l}}$.)
\end{lm}

\begin{defi}\label{boundcl}
\rm Let $(B,\rho,\BBBB)$ be a local contact algebra. A cluster
$\s$ in $(B,C_\rho)$ (see \ref{Alexprn}) is called {\em bounded}\/
if $\s\cap\BBBB\nes$. The set of all bounded clusters in
$(B,C_\rho)$ will be denoted by $\BClu(B,\rho,\BBBB)$. An
ultrafilter $u$ in $B$ is called a {\em bounded ultrafilter}\/ if
$u\cap\BBBB\nes$.
\end{defi}

\begin{nota}\label{compregn}
\rm Let $(X,\tau)$ be a topological space. We will denote by
$CR(X,\tau)$ the family of all compact regular closed subsets of
$(X,\tau)$.  We will often  write  $CR(X)$ instead of
$CR(X,\tau)$.
\end{nota}

\begin{fact}\label{stanlocn}
Let $(X,\tau)$ be a locally compact Hausdorff space. Then the
triple
$$(RC(X,\tau),\rho_{(X,\tau)}, CR(X,\tau))$$
(see \ref{rct} for $\rho_{(X,\tau)}$)
 is
a complete local contact algebra   (\cite{R}). It is called a {\em
standard local contact algebra}.

For every $x\in X$, $\s_x$ is a bounded cluster in
$(RC(X),C_{\rho_X})$ (see (\ref{sxvx}) and (\ref{crho}) for the
notations) (\cite{VDDB}).
\end{fact}

\begin{nist}\label{ladj}
\rm Let $\p:A\lra B$ be an (order-preserving) map between posets,
$A$ has all meets and $\p$ preserves them. Then, by the Adjoint
Functor Theorem (see, e.g., \cite{J}), $\p$ has a left adjoint; it
will be denoted by $\p_\LAM $. Hence $\p_\LAM :B\lra A$ is the
unique order-preserving map such that, for all $a\in A$ and all
$b\in B$, $b\le \p(a)$ iff $\p_\LAM (b)\le a$ (i.e. the pair
$(\pl,\p)$ forms a Galois connection between posets $B$ and $A$).
Equivalently, $\p_\LAM :B\lra A$ is the unique order-preserving
map such that the following two conditions are fulfilled:

\smallskip

\noindent ($\LAM$1) $\fa b\in B$, $\p(\p_\LAM (b))\ge b$;\\
($\LAM$2) $\fa a\in A$, $\p_\LAM (\p(a))\le a$.

\smallskip

\noindent It is well known that $\p\circ\pl\circ\p=\p$,
$\pl\circ\p\circ\pl=\pl$,
\begin{equation}\label{presjoins}
\pl \mbox{ preserves all joins which exist in } B
\end{equation}
and, for all $b\in B$,
\begin{equation}\label{presjoinsf}
\pl(b)=\bigwedge\{a\in A\st\p(a)\ge b\}.
\end{equation}
\noindent Further, $\p$ is an injection iff
\begin{equation}\label{phiin}
\pl(\p(a))=a, \fa a\in A;
\end{equation}
$\p$ is a surjection iff
\begin{equation}\label{phiso}
\p(\pl(b))=b, \fa b\in B.
\end{equation}
\noindent Note that if $\p(0)=0$ then:\\
(a) $\pl(0)=0$ (use ($\LAM$2)), and\\
(b) $\pl(b)\neq 0$, for every $b\in B\stm\{0\}$ (use ($\LAM$1)).

Recall that if $\p\ap:B\lra C$ is a map between posets, $B$ has
all meets and $\p\ap$ preserves them, then
$(\p\ap\circ\p)_\LAM=\pl\circ\pl\ap$.

Finally, if $\psi:A\lra B$ is an (order-preserving) map between
posets, $A$ has all joins and $\psi$ preserves them, then, by the
Adjoint Functor Theorem, $\psi$ has a right adjoint; it will be
denoted by $\psi_P$; $\psi_P:B\lra A$ preserves all meets which
exist in $B$; setting $\p=\psi_P$, we have that $\psi=\pl$.
\end{nist}

\begin{fact}\label{L2rave}
 If $A$ and $B$ are Boolean algebras, $\p:A\lra B$ is a
Boolean homomorphism, $A$ has all meets and $\p$ preserves them,
then:\\
\noindent (a) $\fa a\in A$  and $\fa b\in B$, $\p(a)\we b=0$
implies  $a\we\p_\LAM (b)=0$;\\
\noindent (b) $\fa a\in A$  and  $\fa b\in B$, $\pl(\p(a)\we
b)=a\we\pl(b)$.
\end{fact}

\doc (a) Let $a\in A$, $b\in B$ and $\p(a)\we b=0$. Suppose that
$a\we\pl(b)\neq 0$. Put $c=a\we \pl(b)$. If $\p(c)\we b=0$ then
$b\le\p(c^*)$ and hence $\pl(b)\le c^*$; therefore $c\le c^*$,
i.e. $c=0$,  a contradiction. Thus $\p(c)\we b\neq 0$. This
implies that $\p(a)\we b\neq 0$, a contradiction. Therefore,
$a\we\pl(b)=0$.

\smallskip

\noindent(b) Obviously, $\pl(\p(a)\we b)\le\pl(\p(a))\we\pl(b)\le
a\we\pl(b)$ (by ($\LAM$2) (see \ref{ladj})). Hence, we need only
to show that $\pl(\p(a)\we b)\ge a\we\pl(b)$. By
(\ref{presjoinsf}) (see \ref{ladj}), we have to prove that $a\we
\pl(b)\le\bigwedge\{c\in B\st \p(c)\ge \p(a)\we b\}$. Let $c\in B$
and $\p(c)\ge \p(a)\we b$. We will show that $a\we\pl(b)\le c$.
Using (a) and ($\LAM$1) (see \ref{ladj}), we obtain that:
$a\we\pl(b)\le c\lr c^*\we a\we\pl(b)=0\lr \p(c^*\we a)\we b=0\lr
\p(c)^*\we\p(a)\we b=0\lr \p(a)\we b\le\p(c)$. Thus $a\we\pl(b)\le
c$. Hence (b) is proved.
\sqs

For all undefined here notions and notations see \cite{J, AHS, E,
NW, Si}.


\section{Some New Duality Theorems}

The next  theorem was proved by Roeper \cite{R}.  We will give a
sketch of its proof; it follows the plan of the proof presented in
\cite{VDDB}. The notations and the facts stated here will be used
later on.

\begin{theorem}\label{roeperl}{\rm (P. Roeper \cite{R})}
There exists a bijective correspondence between the class of all
(up to isomorphism) complete local contact algebras and the class
of all (up to homeomorphism) locally compact Hausdorff spaces.
\end{theorem}

\noindent{\em Sketch of the Proof.}~ (A) Let $(X,\tau)$ be a
locally compact Hausdorff space. We put
\begin{equation}\label{psit1}
\Psi^t(X,\tau)=(RC(X,\tau),\rho_{(X,\tau)},CR(X,\tau))
\end{equation}
(see \ref{stanlocn} and \ref{compregn} for the notations).

\noindent(B)~ Let $\unl{B}_{\, l}=(B,\rho,\BBBB)$ be a complete
local contact algebra. Let $C=C_\rho$ be the Alexandroff extension
of $\rho$ (see \ref{Alexprn}). Then, by  \ref{Alexprn}, $(B,C)$ is
a complete normal contact algebra. Put $X=\CL(B,C)$ and let $\TT$
be the topology on $X$ having as a closed base the family
$\{\l_{(B,C)}(a)\st a\in B\}$ where, for every $a\in B$,
\begin{equation}\label{h}
\l_{(B,C)}(a) = \{\s \in X\st  a \in \s\}.
\end{equation}
Sometimes we will write simply $\l_B$ instead of $\l_{(B,C)}$.

\noindent Note that
\begin{equation}\label{intha}
X\stm \l_B(a)= \int(\l_B(a^*)),
\end{equation}
\begin{equation}\label{ee}
\mbox{the family } \{\int(\l_B(a))\st a\in B\} \mbox{ is an open
base of }(X,\TT)
\end{equation}
and, for every $a\in B$,
\begin{equation}\label{haregcl}
\l_B(a)\in RC(X,\TT).
\end{equation}
It can be proved that
\begin{equation}\label{isom}
\l_B:(B,C)\lra (RC(X),\rho_X) \mbox{ is a CA-isomorphism.}
\end{equation}
Further,
\begin{equation}\label{xcomp}
(X,\TT) \mbox{ is a compact Hausdorff space.}
\end{equation}

\noindent(B1)~ Let $1\in\BBBB$. Then $C=\rho$ and $\BBBB=B$, so
that $(B,\rho,\BBBB)=(B,C,B)=(B,C)$ is a complete normal contact
algebra (see \ref{conaln}), and we put
\begin{equation}\label{phiapcn}
\Psi^a(B,\rho,\BBBB)=\Psi^a(B,C,B)=\Psi^a(B,C)=(X,\TT).
\end{equation}

\medskip

\noindent(B2)~ Let $1\not\in\BBBB$. Then, by Lemma \ref{neogrn},
the set $\s_\infty=\{b\in B\st b\not\in\BBBB\}$ is a cluster in
$(B,C)$ and, hence, $\s_\infty\in X$.  Let $L=X\stm\{\s_\infty\}$.
Then
\begin{equation}\label{L}
L=\BClu(B,\rho,\BBBB), \mbox{ i.e. } L \mbox{ is the set of all
bounded clusters of } (B,C_\rho)
\end{equation}
(sometimes we will write $L_{\unl{B}_{\, l}}$ or $L_B$ instead of
$L$);
 let the topology $\tau(=\tau_{\unl{B}_{\, l}})$ on $L$ be the
subspace topology, i.e. $\tau=\TT|_L $. Then $(L,\tau)$ is a
locally compact Hausdorff space. We put
\begin{equation}\label {phiapc}
\Psi^a(B,\rho,\BBBB)=(L,\tau).
\end{equation}

Let
\begin{equation}\label{hapni}
\l^l_{\unl{B}_{\, l}}(a)=\l_{(B,C_\rho)}(a)\cap L,
\end{equation}
for each $a\in B$. We will write simply $\l^l_B$ (or even
$\l_{(A,\rho,\BBBB)}$ when $\BBBB\neq A$) instead of
$\l^l_{\unl{B}_{\, l}}$ when this does not lead to ambiguity. One
can show that:

\smallskip

\noindent (I) $L$ is a dense subset of $X$;\\
(II) $\l^l_B$ is an isomorphism of the Boolean algebra $B$ onto
the
Boolean algebra $RC(L,\tau)$;\\
(III) $b\in\BBBB$ iff $\l^l_B(b)\in CR(L)$;\\
(IV) $a\rho b$ iff $\l^l_B(a)\cap \l^l_B(b)\neq\ems$.\\
Hence, $X$ is the Alexandroff (i.e. one-point) compactification of
$L$ and
\begin{equation}\label{hapisom}
\l^l_B: (B,\rho,\BBBB)\lra (RC(L),\rho_L, CR(L)) \mbox{ is an
LCA-isomorphism.}
\end{equation}
 Note also that for every $b\in B$,
\begin{equation}\label{intl}
\int_{L_B}(\l^l_B(b))=L_B\cap\int_X(\l_B(b)).
\end{equation}

\medskip

\noindent(C)~ For every CLCA $(B,\rho,\BBBB)$ and every $a\in B$,
set
\begin{equation}\label{lbg}
\l^g_{\unl{B}_{\,l}}(a)=\l_{(B,C_\rho)}(a)\cap\Psi^a(B,\rho,\BBBB).
\end{equation}
We will write simply $\l^g_B$ instead of $\l^g_{\unl{B}_{\, l}}$
when this does not lead to ambiguity.
 Thus, when $1\in\BBBB$,
we have that $\l^g_B=\l_B$, and  if $1\nin\BBBB$  then
$\l^g_B=\l^l_B$. Hence, by (\ref{isom}) and (\ref{hapisom}), we
get that
\begin{equation}\label{hapisomn}
\l^g_B: (B,\rho,\BBBB)\lra (\Psi^t\circ\Psi^a)(B,\rho,\BBBB)
\mbox{ is an LCA-isomorphism.}
\end{equation}

With the next assertion we specify (\ref{ee}):
\begin{equation}\label{eel}
\mbox{the family } \{\int_{\Psi^a(B,\rho,\BBBB)}(\l_B^g(a))\st
a\in \BBBB\} \mbox{ is an open base of } \Psi^a(B,\rho,\BBBB).
\end{equation}

\noindent(D)~~ Let $(X,\tau)$ be a compact Hausdorff space. Then
it can be proved that the map
\begin{equation}\label{nison}
t_{(X,\tau)}:(X,\tau)\lra\Psi^a(\Psi^t(X,\tau)),
\end{equation}
defined by  $t_{(X,\tau)}(x)=\{F\in RC(X,\tau)\st x\in
F\}(=\s_x)$, for all $x\in X$, is a homeomorphism (we will also
write simply $t_X$ instead of $t_{(X,\tau)}$).

Let  $(L,\tau)$ be a non-compact locally compact Hausdorff space.
Put $B=RC(L,\tau)$, $\BBBB=CR(L,\tau)$ and $\rho=\rho_L$. Then
$(B,\rho,\BBBB)=\Psi^t(L,\tau)$ and $1\nin\BBBB$ (here $1=L$). It
can be shown that the map
\begin{equation}\label{homeo}
t_{(L,\tau)}:(L,\tau)\lra\Psi^a(\Psi^t(L,\tau)),
\end{equation}
defined by  $t_{(L,\tau)}(x)=\{F\in RC(L,\tau)\st x\in
F\}(=\s_x)$, for all $x\in L$, is a homeomorphism.

Therefore $\Psi^a(\Psi^t(L,\tau))$ is homeomorphic to $(L,\tau)$
and $\Psi^t(\Psi^a(B,\rho,\BBBB))$ is LCA-isomorphic to
$(B,\rho,\BBBB)$.
 \sqs

\begin{cor}\label{vriesb}{\rm (De Vries \cite{dV})}
There exists a bijective correspondence between the class of all
(up to isomorphism) complete normal contact algebras and the class
of all (up to homeomorphism)  compact Hausdorff spaces.
\end{cor}

\doc The restriction of the correspondence $\Psi^a$, defined in the
proof of Theorem \ref{roeperl},
 to the class
of all complete normal contact algebras generates the required
bijective correspondence (see (B1) in the proof of \ref{roeperl}).
\sqs

\begin{defi}\label{semiopen}{\rm (\cite{MP})}
\rm
 A continuous map $f:X\lra Y$ is called {\em quasi-open\/} if for
every non-empty open subset $U$ of $X$, $\int(f(U))\nes$ holds.
\end{defi}

Every closed irreducible map $f:X\lra Y$ is quasi-open (because,
for every non-empty open subset $U$ of $X$, $f^{\rm \#}(U)(=\{y\in
Y\st f\inv(y)\sbe U\})$ is a non-empty open subset of $Y$
(\cite{P})).

Recall that a function $f:X\lra Y$ is called  {\em skeletal}\/
(\cite{MR}) if
\begin{equation}\label{ske}
\int(f\inv(\cl (V)))\sbe\cl(f\inv(V))
\end{equation}
for every open subset  $V$  of $Y$. As it is noted in \cite{MR}, a
continuous map $f:X\lra Y$ is skeletal iff $f\inv(\Fr(V))$ is
nowhere dense in $X$, for every open subset $V$ of $Y$. Clearly, a
function $f:X\lra Y$ is  skeletal  iff $\int(f\inv(\Fr(V)))=\ems$,
for every open subset $V$ of $Y$. The next assertion can be easily
proved:

\begin{lm}\label{skel}
A function $f:X\lra Y$ is  skeletal iff\/ $\int(\cl(f(U)))\nes$,
for every  non-empty  open subset $U$ of $X$.
\end{lm}

\doc ($\Rightarrow$) Let $U\/$ be a  non-empty open subset of $X$.
Suppose that $\int(\cl(f(U)))=\ems$. Set $V=Y\stm \cl(f(U))$. Then
$\Fr(V)=Y\stm V=\cl(f(U))$ and hence $U\sbe f\inv(\Fr(V))$. Thus
$\int(f\inv(\Fr(V)))\nes$, a contradiction. Therefore,
$\int(\cl(f(U)))\nes$.

\smallskip

($\Leftarrow$) Let $V$ be an open subset of $Y$. Suppose that
$U=\int(f\inv(\Fr(V)))$ is a non-empty set. Then $\ems\neq
\int(\cl(f(U)))\sbe\Fr(V)=\cl(V)\stm V$, which is impossible.
Hence $\int(f\inv(\Fr(V)))=\ems$. So, $f$ is a skeletal map.
\sqs

A topological space $(X,\tau)$ is said to be $\pi$-{\em regular}\/
 if  for each non-empty $U\in\tau$ there exists
   a non-empty $V\in\tau$ such that $\cl(V)\sbe U$. The semiregular
   $\pi$-regular spaces are exactly the  {\em weakly
regular}\/ spaces of D\"{u}ntsch and Winter (\cite{DW}).

\begin{cor}\label{semsk}
(a) Every quasi-open map is skeletal.

\noindent (b) Let $X$ be a $\pi$-regular space and $f:X\lra Y$ be
a closed map. Then $f$ is quasi-open iff $f$ is skeletal.
\end{cor}

\doc (a) It follows from \ref{skel}.

\smallskip

\noindent (b) Let $f$ be skeletal and closed. Take an open
non-empty subset $U$ of $X$. Then there exists an open non-empty
subset $V$ of $X$ such that $\cl(V)\sbe U$. Using \ref{skel}, we
obtain that $\int(f(U))\spe\int(f(\cl(V)))=\int(\cl(f(V)))\nes$.
Therefore, $f$ is a quasi-open map.
\sqs

\begin{lm}\label{mrro}
Let $f:X\lra Y$ be a continuous map. Then the
following conditions are equivalent:\\
(a) $f$ is a skeletal map;\\
(b) For every $F\in RC(X)$, $\cl(f(F))\in RC(Y)$.
\end{lm}

\doc (a)$\Rightarrow$(b) Let $f$ be a skeletal map, $F\in
RC(X)$ and $F\nes$. Set $U=\int(F)$. Then $U\nes$. Hence, by
\ref{skel}, $V=\int(\cl(f(U)))\nes$. We will show that
\begin{equation}\label{clff}
\cl(f(F))=\cl(V).
\end{equation}
Note that, by the continuity of $f$, $\cl(f(F))=\cl(f(U))$. Now
suppose that $f(U)\not\sbe\cl(V)$. Then there exists $y\in
f(U)\stm \cl(V)$. Hence there exists an open neighborhood $O_1$ of
$y$ in $Y$ such that $O_1\cap V=\ems$. Thus $\cl(O_1)\cap V=\ems$.
There exists $x\in U$ such that $y=f(x)$. Since $f$ is continuous,
there exists an open neighborhood $O$ of $x$ in $X$  such that
$x\in O\sbe U$ and $f(O)\sbe O_1$. Then $\cl(f(O))\sbe\cl(O_1)$
and thus $\cl(f(O))\cap V=\ems$. Since, by \ref{skel},
$\ems\neq\int(\cl(f(O)))\sbe
\cl(f(O))\cap\int(\cl(f(U)))=\cl(f(O))\cap V=\ems$, we obtain a
contradiction. Therefore $f(U)\sbe\cl(V)$ and hence
$\cl(f(U))\sbe\cl(V)$. Since the converse inclusion is obvious,
(\ref{clff}) is established. Thus, $\cl(f(F))\in RC(Y)$.

\medskip

\noindent (b)$\Rightarrow$(a)  Let $U$ be a non-empty open subset
of $X$. Then $F=\cl(U)\in RC(X)$. Hence $\cl(f(F))\in RC(Y)$.
Since $F\nes$, we obtain that $\int(\cl(f(F)))\nes$. Now, using
the continuity of $f$, we get that $\int(\cl(f(U)))\nes$.
Therefore, by \ref{skel}, $f$ is a skeletal map.
\sqs

The next lemma generalizes the well-known result of Ponomarev
\cite{P} that the regular closed sets are preserved by the closed
irreducible maps.

\begin{lm}\label{pon}
Let $f:X\lra Y$ be a closed map and $X$ be a $\pi$-regular space.
Then the
following conditions are equivalent:\\
(a) $f$ is a quasi-open map;\\
(b) For every $F\in RC(X)$, $f(F)\in RC(Y)$.
\end{lm}

\doc (a)$\Rightarrow$(b) It follows from \ref{semsk}(a) and
\ref{mrro}.

\medskip

\noindent (b)$\Rightarrow$(a) It follows from \ref{semsk}(b) and
\ref{mrro}. Note that the $\pi$-regularity of $X$ is used only in
the proof of this implication.
\sqs

\begin{cor}\label{ponco}
If $f:X\lra Y$ is a quasi-open closed map then $f(X)\in RC(Y)$.
\end{cor}

\begin{rems}\label{remhj}
\rm
In \cite{HJ}, Henriksen and Jerison regarded  functions $f:X\lra
Y$ for which
\begin{equation}\label{skehj}
\cl(\int(f\inv(F)))=\cl(f\inv(\int(F))) \mbox{ for every } F\in
RC(Y).
\end{equation}
Clearly, every continuous skeletal map $f:X\lra Y$ satisfies
(\ref{skehj}) (\cite{MR}). Hence, by \ref{semsk}(a), every
quasi-open map $f:X\lra Y$ satisfies (\ref{skehj}) (\cite{PS}).

Functions $f:X\lra Y$ (not necessarily continuous) satisfying
condition (\ref{ske}) for every $V\in RO(X)$ are called {\em
HJ-maps}\/ in \cite{MR}. Obviously, every continuous HJ-map
$f:X\lra Y$ satisfies (\ref{skehj}). As it is noted in \cite{MR},
the composition of two continuous HJ-maps needs not
be an
HJ-map, while the composition of two continuous skeletal maps is a
skeletal map. It is clear that the composition of two quasi-open
maps is a quasi-open map.
\end{rems}

\begin{defi}\label{lcatnk}
\rm Let $\SKLC$ be the category of all locally compact Hausdorff
spaces and all continuous skeletal maps between them.

Let $\SKAL$ be the category whose objects are all complete local
contact algebras
 and whose morphisms   $\p:(A,\rho,\BBBB)\lra
(B,\eta,\BBBB\ap)$   are all
  complete Boolean homomorphisms $\p:A\lra B$  satisfying the
following conditions:

\smallskip

\noindent (L1) $\fa a,b\in A$, $\p(a)\eta\p(b)$ implies
$a\rho b$;\\
(L2) $b\in\BBBB\ap$ implies $\pl(b)\in\BBBB$ (see \ref{ladj} for
$\pl$).

\smallskip

It is easy to see that in this  way we have defined categories.

Let us  note  that (L1) is equivalent to the following condition:

\smallskip

\noindent (EL1) $\fa a,b\in B$, $a\eta b$ implies
$\pl(a)\rho\pl(b)$.
\end{defi}


\begin{theorem}\label{maintheoremnk}
The categories $\SKLC$ and $\SKAL$ are dually equivalent.
\end{theorem}

\doc We will define two
contravariant functors $\Psi^a:\SKAL\lra\SKLC$ and
$\Psi^t:\SKLC\lra\SKAL$. On the objects they coincide with the
correspondences $\Psi^a$ and $\Psi^t$, respectively (see
(\ref{psit1}), (\ref{phiapcn}) and (\ref{phiapc}) for them). We
will define $\Psi^a$ and $\Psi^t$ on the morphisms of the
corresponding categories.

Let $f\in\SKLC((X,\tau),(Y,\tau\ap))$. Define
$\Psi^t(f):\Psi^t(Y,\tau\ap)\lra\Psi^t(X,\tau)$ by
\begin{equation}\label{defpsiapnk}
\Psi^t(f)(F)=\cl(f\inv(\int(F))), \fa F\in\Psi^t(Y,\tau\ap).
\end{equation}
Then, by \ref{remhj},
\begin{equation}\label{defpsin12k}
\Psi^t(f)(F)=\cl(\int(f\inv(F))), \fa F\in \Psi^t(Y,\tau\ap).
\end{equation}
Put $\p=\Psi^t(f)$. We will first show that $\p$ is a complete
Boolean homomorphism. Let $\GA$ be a set and $\{F_\g\st
\g\in\GA\}\sbe RC(Y)$. Put $F=\cl(\bigcup\{F_\g\st \g\in\GA\})$.
(Note that $F=\cl(\bigcup\{\int(F_\g)\st\g\in\GA\})$.) Then $F\in
RC(Y)$ and $\bigvee\{F_\g\st \g\in\GA\}=F$. Since $\p$ is an
order-preserving map, we get that
$\p(F)\ge
\bigvee\{\p(F_\g)\st \g\in\GA\}$. We will now prove the converse
inequality. We have that $\p(F)=\cl(f\inv(\int(F)))$. Let $x\in
f\inv(\int(F))$. Then $f(x)\in\int(F)$. Hence, there exist  open
neighborhoods $O$ and $O\ap$ of $f(x)$ in $Y$ such that
$\cl(O\ap)\sbe O\sbe F$. Since $f$ is continuous, there exists an
open neighborhood $U$ of $x$ in $X$ such that $f(U)\sbe O\ap$.
Suppose that there exists an open neighborhood $V$ of $x$ in $X$
such that, for every $\g\in\GA$,
$V\cap\cl(\int(f\inv(F_\g)))=\ems$. Obviously, we can suppose that
$V\sbe U$. Since $f$ is continuous and skeletal, we get, using
\ref{remhj} and (\ref{skehj}), that $V\cap
f\inv(\int(F_\g))=\ems$, for every $\g\in\GA$. Thus,
$f(V)\cap\bigcup\{\int(F_\g)\st\g\in\GA\}=\ems$. Put
$W=\bigcup\{\int(F_\g)\st\g\in\GA\}$. Then $\cl(f(V))\cap W=\ems$
and $\cl(f(V))\sbe \cl(f(U)\sbe\cl(O\ap)\sbe O\sbe F=\cl(W)$. Thus
$\cl(f(V))\sbe\cl(W)\stm W=\Fr(W)$. Since $f$ is skeletal,
\ref{skel} implies that $\int(\cl(f(V)))\nes$ and this leads to a
contradiction. Therefore,
$x\in\cl(\bigcup\{\cl(\int(f\inv(F_\g)))\st\g\in\GA\})$. We have
proved that $\p(F)\sbe\bigvee\{\p(F_\g)\st \g\in\GA\}$. So,
$\p(\bigvee\{F_\g\st \g\in\GA\})=\bigvee\{\p(F_\g)\st \g\in\GA\}$.

Let $F\in RC(Y)$. Then, by (\ref{defpsiapnk}) and
(\ref{defpsin12k}),
$(\p(F))^*=(\cl(f\inv(\int(F))))^*=(\cl(f\inv(Y\stm
F^*)))^*=(\cl(X\stm f\inv(F^*)))^*=\cl(X\stm\cl(X\stm
f\inv(F^*)))=\cl(\int(f\inv(F^*)))$. So, $\p(F^*)=(\p(F))^*$.
Since, obviously, $\p$ preserves the zero and the unit elements,
$\p$ is a complete Boolean homomorphism.

Further, using \ref{mrro}, we can define a map
\begin{equation}\label{psilnk}
\psi:  \Psi^t(X,\tau)\lra  \Psi^t(Y,\tau\ap) \mbox{ by
}\psi(G)=\cl(f(G)), \mbox{ for every } G\in  \Psi^t(X,\tau).
\end{equation}
 Obviously, $\psi$ is an order-preserving map. Since $f$ is a
 continuous
 map, we have that
 for every $F\in RC(Y)$,
 $\psi(\p(F))=\cl(f(\cl(f\inv(\int(F)))))=\cl(f(f\inv(\int(F))))\sbe\cl(\int(F))=F$,
  and, similarly, for every
$G\in RC(X)$,
$\p(\psi(G))=\p(\cl(f(G)))=\cl(\int(f\inv(\cl(f(G)))))\spe
\cl(\int(G))=G$. Hence $\psi$ is a left adjoint to $\p$ (see
\ref{ladj}), i.e.
\begin{equation}\label{leftadjnk}
\psi=\p_\LAM.
\end{equation}

We have to show that $\p$ satisfies conditions (L1) and (L2).
Using (\ref{leftadjnk}), we obtain immediately that (EL1) (and
hence (L1)) and (L2) are fulfilled.

Hence, $\Psi^t(f)$ is a morphism of the category $\SKAL$.

It is obvious that $\Psi^t(id)=id$. Let us show that
$\Psi^t(g\circ f)= \Psi^t(f)\circ\Psi^t(g)$. Indeed, using
continuity of $f$ and $g$, \ref{remhj}  and (\ref{skehj}), we
obtain that $(\Psi^t(f)\circ\Psi^t(g))(F)=
\cl(\int(f\inv(\cl(g\inv(\int(F))))))\spe
\cl(\int(\cl(f\inv(g\inv(\int(F))))))=\Psi^t(g\circ f)(F)$ and
also $(\Psi^t(f)\circ
\Psi^t(g))(F)=\cl(f\inv(\int(\Psi^t(g)(F))))\sbe\cl(\int(f\inv(\cl(\int(g\inv(F))))))
\sbe\cl(\int(f\inv(g\inv(F))))=\Psi^t(g\circ f)(F)$. Thus,
$$  \Psi^t:\SKLC\lra \SKAL$$
is a  contravariant functor.

\smallskip

Let $\p\in\SKAL((A,\rho,\BBBB),(B,\eta,\BBBB\ap))$. Since
$\p:A\lra B$ is a complete Boolean homomorphism, $\p$ has a left
adjoint $\p_\LAM :B\lra A$ (see \ref{ladj}). Set $C=C_\rho$ and
$C\ap=C_\eta$ (see \ref{Alexprn} for the notations). We will write
$``\ll$" for $``\ll_C$" and $``\ll\ap$" for $``\ll_{C\ap}$".

Define now
\begin{equation}\label{psiphink}
\Psi^a(\p):\Psi^a(B,\eta,\BBBB\ap)\lra \Psi^a(A,\rho,\BBBB)
\end{equation}
by the formula
\begin{equation}\label{formnk}
\Psi^a(\p)(\s_u)=\s_{\p\inv(u)},
\end{equation}
where $u\in \Ult(B)$, $\s_u$ is a cluster in $(B,C\ap)$,
$\s_u\cap\BBBB\ap\nes$ and $\s_{\p\inv(u)}$ is a cluster in
$(A,C)$ (see (\ref{sigmau}) and \ref{uniqult} for $\s_u$ and
$\s_{\p\inv(u)}$, and note that, by \ref{conclustth}, any cluster
$\s$ in $(B,C\ap)$ can be written in the form $\s_u$ for some
$u\in \Ult(B)$).

We have to show that $\Psi^a(\p)$ is well defined. Set
$f=\Psi^a(\p)$, $X=\Psi^a(A,\rho,\BBBB)$ and
$Y=\Psi^a(B,\eta,\BBBB\ap)$. Then $X$ is the set of all bounded
clusters of $(A,C)$ and $Y$ is the set of all bounded clusters of
$(B,C\ap)$ (see \ref{boundcl}, (\ref{phiapcn}) and (\ref{L})).

Let us start with the following observation:
\begin{equation}\label{ultbasis}
\mbox{if } u\in\Ult(B) \mbox{ then }\p\inv(u)\in \Ult(A) \mbox{
and } \pl(u) \mbox{ is a basis of } \p\inv(u).
\end{equation}
So, let $u\in\Ult(B)$. Then, obviously, $\p\inv(u)\in \Ult(A)$.
Let us show that $\pl(u)\sbe\p\inv(u)$. Let $b\in u$. Then, by
($\LAM$1), $\p(\pl(b))\ge b$. Hence $\p(\pl(b)\in u$, i.e.
$\pl(b)\in\p\inv(u)$. Therefore, $\pl(u)\sbe\p\inv(u)$. Further,
suppose that there exists $a\in\p\inv(u)$ such that $\pl(b)\not\le
a$ for all $b\in u$. Then $\pl(b)\we a^*\neq 0$ for every $b\in
u$. Hence, by \ref{L2rave}(a), $b\we\p(a^*)\neq 0$ for every $b\in
u$. Since $u\in\Ult(B)$, we obtain that $\p(a^*)\in u$. Thus both
$\p(a)$ and $\p(a)^*$ are elements of $u$,  a contradiction.
Therefore, $\pl(u)$ is a basis of the ultrafilter $\p\inv(u)$.

Obviously, (\ref{ultbasis}) implies that
\begin{equation}\label{plu}
\fa u\in\Ult(B), \s_{\p\inv(u)}=\s_{\pl(u)},
\end{equation}
where $\s_{\p\inv(u)}$ and $\s_{\pl(u)}$ are clusters in $(A,C)$
(see \ref{uniqult} for the notations).

Let $\s$ be a cluster in $(B,C\ap)$. Then the following holds:
\begin{equation}\label{bstar}
\mbox{if } \s\cap\BBBB\ap\nes\mbox{ then there exists } b\in
\BBBB\ap \mbox{ such that } b^*\nin\s.
\end{equation}
Indeed, let $b_0\in\s\cap\BBBB\ap$. Since $b_0\ll_\eta 1$, (BC1)
implies that there exists $b\in\BBBB\ap$ such that $b_0\ll_\eta
b$. Then $b_0(-\eta)b^*$ and since $b_0\in\BBBB\ap$, we obtain
that $b_0(-C\ap)b^*$. Thus $b^*\nin\s$.

Let us now show  that
\begin{equation}\label{boundult}
\mbox{if } u\in\Ult(B) \mbox{ and } \s_u\cap\BBBB\ap\nes \mbox{
then } u\cap\BBBB\ap\nes
\end{equation}
(here, of course, $\s_u$ is a cluster in $(B,C\ap)$). Indeed,  by
(\ref{bstar}), there exists $a\in\BBBB\ap$ such that
$a^*\nin\s_u$. Hence $a\in u\cap\BBBB\ap$. So, (\ref{boundult}) is
proved.

Let $u,v\in \Ult(B)$, $\s_u=\s_v$ and $\s=\s_u(=\s_v)$ be bounded.
We will prove that $\s_{\p\inv(u)}=\s_{\p\inv(v)}$. Indeed, by
(\ref{boundult}), there exists $c\in u\cap \BBBB\ap$. Let $a\in u$
and $b\in v$. Then $a\we c\in u\cap\BBBB\ap$ and $(a\we c)C\ap b$.
Thus $(a\we c)\eta b$. Hence, by (EL1), $\pl(a\we c)\rho\pl(b)$.
Therefore, $\pl(a\we c)C\pl(b)$. Thus $\pl(a)C\pl(b)$. Since this
is true for every $a\in u$ and every $b\in v$, we obtain, using
(\ref{ultbasis}) and (\ref{sigmau}), that $\pl(u)\sbe\s_{\pl(v)}$.
Then, by \ref{uniqult} and (\ref{ultbasis}),
$\s_{\pl(u)}=\s_{\pl(v)}$. Using (\ref{plu}), we get that
$\s_{\p\inv(u)}=\s_{\p\inv(v)}$.

Now, using (\ref{plu}), we obtain  that
\begin{equation}\label{philbn}
\mbox{if } \s\in Y \mbox{ and } b\in\s \mbox{ then }\p_\LAM(b)\in
f(\s).
\end{equation}
Indeed, by \ref{conclustth}, there exists $u\in \Ult(B)$ such that
$b\in u\sbe\s$, and hence $\s=\s_u$. Thus, by (\ref{plu}),
$f(\s)=\s_{\pl(u)}$.  Therefore $\p_\LAM(b)\in f(\s)$. So,
(\ref{philbn}) is proved.

Let us show that for every $\s\in \CL(B,C\ap)$,
\begin{equation}\label{sbounn}
\s\cap\BBBB\ap\nes\mbox{ implies that } f(\s)\cap\BBBB\nes.
\end{equation}
Indeed, let $\s\in \CL(B,C\ap)$ and $b\in\s\cap\BBBB\ap$. Then, by
(\ref{philbn}), $\p_\LAM(b)\in f(\s)$. Since, by (L2),
$\p_\LAM(b)\in\BBBB$, we obtain that $f(\s)\cap\BBBB\nes$.

So, the function $f$ is well defined on $Y$ and $f(Y)\sbe X$. We
have to show that $f$ is continuous and skeletal.

Note first that, using (\ref{intha}) and (\ref{intl}), we get
readily that for every $a\in A$,
\begin{equation}\label{inthal}
X\stm \l_A^g(a)= \int_X(\l_A^g(a^*)).
\end{equation}

Further, using (\ref{hapisomn}) and \ref{rct}, one can easily show
that for all $a,b\in A$,
\begin{equation}\label{abipn}
a\ll_\rho b \mbox{ implies that } \l_A^g(a) \sbe\int_X(\l_A^g(b)).
\end{equation}

Note also that if $\s$ is a cluster in $(B,C\ap)$ then
\begin{equation}\label{b1b2}
b_1^*,b_2^*\nin\s\mbox{ implies that } b_1\we b_2\in\s\mbox{ and }
(b_1\we b_2)^*\nin\s.
\end{equation}
Indeed, if $b_1^*,b_2^*\nin\s$ then, by (K2), $b_1^*\vee
b_2^*\nin\s$, i.e. $(b_1\we b_2)^*\nin\s$;  hence $b_1\we
b_2\in\s$.

Let us now prove that for every $b\in \BBBB\ap$,
\begin{equation}\label{raven2n}
f(\l_B^g(b))=\l_A^g(\p_\LAM (b))
\end{equation}
(note that  $b\in\BBBB\ap$ implies that $\l_B(b)\sbe Y$ and
$\pl(b)\in\BBBB$ (by (L2));  thus we have also that $\l_A(\p_\LAM
(b))\sbe X$; hence (\ref{raven2n}) can be written as
$f(\l_B(b))=\l_A(\p_\LAM (b))$).   Since $\p(0)=0$, we have, by
\ref{ladj}, that $\pl(0)=0$ and $\pl(b)\not=0$ for any $b\not=0$.
Hence, (\ref{raven2n}) is true for $b=0$.

 Let $b\in \BBBB\ap\stm\{0\}$ and $\s\in f(\l_B(b))$. Then there
exists $\s\ap\in \l_B(b)$ such that $f(\s\ap)=\s$. Hence
$b\in\s\ap$ and thus, by (\ref{philbn}), $\p_\LAM (b)\in
f(\s\ap)=\s$. Therefore we get that $\s\in \l_A(\p_\LAM (b))$. So,
$f(\l_B(b))\sbe \l_A(\p_\LAM (b))$. Conversely, let $b\in
\BBBB\ap\stm\{0\}$ and $\s\in \l_A(\p_\LAM (b))$, i.e. $\p_\LAM
(b)\in\s$. Then, by \ref{conclustth}, there exists $u\in \Ult(A)$
such that $\p_\LAM (b)\in u\sbe\s$, and hence, by \ref{uniqult},
$\s=\s_u$. Let us show that $\p(u)\cup\{b\}$  has the finite
intersection property. Since $\p(u)$ is closed under finite meets,
it is enough to prove that $b\we\p(a)\neq 0, \fa a\in u$. Indeed,
suppose that there exists $a_0\in u$ such that $b\we \p(a_0)=0$.
Then, by \ref{L2rave}(a), we will have that $\p_\LAM (b)\we
a_0=0$. This is, however, impossible, since $\p_\LAM (b)\in u$.
So, there exists an ultrafilter $v$ in $B$ such that
$v\spe\p(u)\cup\{b\}$. Set $\s\ap=\s_v$. Then $\s\ap$ is a cluster
in $(B,C\ap)$ (see \ref{conclustth}) and since $v\sbe\s\ap$, we
have that $b\in\s\ap$. Hence $\s\ap\in \l_B(b)$. Further,
$f(\s\ap)=\s$. Indeed, since $\p(u)\sbe v$, we have that
$u\sbe\p\inv(v)$; thus $u=\p\inv(v)$ and hence
$\s=\s_u=\s_{\p\inv(v)}=f(\s_v)=f(\s\ap$). Therefore
$\s=f(\s\ap)\in f(\l_B(b))$. So, (\ref{raven2n}) is proved.

We are ready to show that $f$ is a continuous function.

Let $\s\in Y$, $\s\ap=f(\s)$, $a\in A$ and
$\s\ap\in\int_X(\l_A^g(a))$ (we use (\ref{eel})). Then, by
(\ref{inthal}), $a^*\nin\s\ap$. By \ref{cluendcor}, there exists
$a_1\in A$ such that $a^*\ll a_1^*$ and $a_1^*\nin\s\ap$. Then
$a_1\in v$, for every $v\in\Ult(A)$ such that $v\sbe\s\ap$. Thus,
using (\ref{ultbasis}), we obtain that for every $u\in\Ult(B)$
such that $u\sbe\s$, there exists $b_u\in u$ with $\pl(b_u)\le
a_1$. Set $b=\bigvee\{b_u\st u\in\Ult(B), u\sbe\s\}$. Then, by
\ref{ladj}, $\pl(b)=\bigvee\{\pl(b_u)\st u\in\Ult(B), u\sbe\s\}$.
Hence $\pl(b)\le a_1$. Suppose that $b^*\in\s$. Then
\ref{conclustth} implies that there exists $u\in\Ult(B)$ such that
$b^*\in u\sbe\s$. Since $b\in u$ (because $b_u\in u$ and $b_u\le
b$), we obtain a contradiction. Hence $b^*\nin\s$.  Since $\s$ is
a bounded cluster, (\ref{bstar}) implies that there exists
$c\in\BBBB\ap$ such that $c^*\nin\s$. Set $d=b\we c$. Then
$d\in\BBBB\ap$ and $d^*\nin\s$ (by (\ref{b1b2})). Now, using (L2),
(\ref{abipn}) and (\ref{raven2n}), we obtain that
$f(\int_Y(\l_B^g(d)))\sbe
f(\l_B^g(d))=\l_A^g(\pl(d))\sbe\l_A^g(\pl(b))\sbe\l_A^g(a_1)\sbe
\int_X(\l_A^g(a))$. Since $\s\in\int_Y(\l_B^g(d))$, we get that
\begin{equation}\label{fcontk}
f:Y\lra X \mbox{ is a continuous function.}
\end{equation}

We will now show that $f$ is a skeletal map. Since $f$ is
continuous, it is enough to prove, by \ref{skel}, that
$\int_X(f(\cl(U)))\nes$ for every non-empty open subset $U$ of
$Y$. Hence, by (\ref{eel}) and (\ref{hapisomn}), we have to show
that $\int_X(f(\l_B^g(b))\nes$, for every $b\in
\BBBB\ap\stm\{0\}$.

Suppose that there exists $b\in \BBBB\ap\stm\{0\}$ such that
$\int_X(f(\l_B^g(b)))=\ems$. Then $X\stm f(\l_B^g(b))$ is dense in
$X$. Using (\ref{raven2n}), we obtain that $X\stm\l_A^g(\pl(b))$
is dense in $X$. Thus, by (\ref{inthal}),
$\int(\l_A^g((\pl(b))^*))$ is dense in $X$. Now, (\ref{hapisomn})
implies that $\l_A^g((\pl(b))^*)=X$. Therefore, by
(\ref{hapisomn}), $(\pl(b))^*=1$. Then $\pl(b)=0$ and hence $b=0$
(by \ref{ladj}),  a contradiction. Hence,
\begin{equation}\label{skm}
f:Y\lra X \mbox{ is a skeletal map.}
\end{equation}

So, we have proved that
$\Psi^a(\p)\in\SKLC(\Psi^a(B,\eta,\BBBB\ap),
\Psi^a(A,\rho,\BBBB))$. Thus $\Psi^a$ is well defined on the
morphisms of the category $\SKAL$.

It is easy to see that $\Psi^a$ preserves the identity maps and
that $ \Psi^a(\p_1\circ\p_2)=  \Psi^a(\p_2)\circ  \Psi^a(\p_1)$.
Thus,
$$  \Psi^a:\SKAL\lra \SKLC$$
is a  contravariant functor.

We will prove that $\Psi^a\circ  \Psi^t\cong Id_{\SKLC}$ (where
$``\cong''$ means $``$naturally isomorphic" and $Id$ is the
identity functor).

We will show that
\begin{equation}\label{alphani1nk}
t:Id_{\SKLC}\lra   \Psi^a\circ  \Psi^t,
\end{equation}
 defined by
\begin{equation}\label{alphanink}
 t(X,\tau)=t_{(X,\tau)}, \fa (X,\tau)\in\card\SKLC ,
\end{equation}
is the required natural isomorphism (see (\ref{nison}) and
(\ref{homeo}) for the definition of $t_{(X,\tau)}$).

Let $f\in\SKLC((X,\tau),(Y,\tau\ap))$ and $\fs=
\Psi^a(\Psi^t(f))$. We have to show that $\fs\circ t_X=t_Y\circ
f$. Let $x\in X$. Then $\fs(t_X(x))=\fs(\s_x)$ and $(t_Y\circ
f)(x)=\s_{f(x)}$. We will prove that
\begin{equation}\label{fsxnk}
\fs(\s_x)=\s_{f(x)}.
\end{equation}
Note first that
\begin{equation}\label{usxnk}
\mbox{if } u\in \Ult(RC(X)),\ x\in X \mbox{ and } u\supset \nu_x
\mbox{ then } u\subset\s_x
\end{equation}
(see (\ref{sxvx}) for $\nu_x$). Indeed, let $F\in u$ and suppose
that $x\nin F$. Then $x\in X\stm F=\int(F^*)$ and hence $F^*\in
\nu_x$. Thus $F^*\in u$,  a contradiction. So, $u\subset\s_x$.

Set $\Psi^t(f)=\p$. Let $x\in X$. Since $\nu_x$ is a filter in
$RC(X)$ (see \ref{sxcluster}), there exists $u\in \Ult(RC(X))$
such that $\nu_x\sbe u$. Then, by (\ref{usxnk}) and \ref{uniqult},
$\s_x=\s_u$. Hence $\fs(\s_x)=\s_{\p\inv(u)}$. We will now show
that $\nu_{f(x)}\sbe\p\inv(u)$. Indeed, let $G\in RC(Y)$ and
$f(x)\in\int_Y(G)$. Then, by the continuity of $f$, $x\in
f\inv(\int_Y(G))\sbe\int_X(\p(G))$ (see (\ref{defpsin12k})). Thus
$\p(G)\in \nu_x\sbe u$. Hence $G\in\p\inv(u)$. Therefore
$\nu_{f(x)}\sbe \p\inv(u)$.
 Then, by
(\ref{usxnk}) and \ref{uniqult},
$\s_{f(x)}=\s_{\p\inv(u)}=\fs(\s_x)$. So, we have proved that
$\fs(t_X(x))=t_Y(f(x))$, for every $x\in X$. Hence,
 $t$ is a natural isomorphism.

Finally, we will prove that $\Psi^t\circ  \Psi^a\cong Id_{\SKAL}$.

We will show that
\begin{equation}\label{be1nk}
\l: Id_{\SKAL}\lra  \Psi^t\circ  \Psi^a,\mbox{ where }
\l(A,\rho,\BBBB)=\l_A^g, \fa (A,\rho,\BBBB)\in\card\SKAL
\end{equation}
(see (\ref{lbg}) for $\l_A^g$), is the required natural
isomorphism.

Let $(A,\rho,\BBBB)\in\card\SKAL$. Using (\ref{hapisomn}), it is
easy to see that
\begin{equation}\label{hbmorphink}
\l_A^g:(A,\rho,\BBBB)\lra  \Psi^t(  \Psi^a(A,\rho,\BBBB)) \mbox{
is an } \SKAL\mbox{-isomorphism.}
\end{equation}

Let $\p\in\SKAL((A,\rho,\BBBB),(B,\eta,\BBBB\ap))$ and
$\ps=\Psi^t(\Psi^a(\p))$. We have to prove that
$\l_B^g\circ\p=\ps\circ \l_A^g$. Set $f= \Psi^a(\p)$. Let $a\in
A\stm\{0\}$. Put $F=\l_A^g(a)$ and $G=\l_B^g(\p(a))$. We have to
show that $\ps(F)=G$, i.e. that
$G=\cl(f\inv(\int(F)))(=\cl(\int(f\inv(F))))$. Let $\s\in G$. Then
$\p(a)\in\s$ and $\s\cap\BBBB\ap\nes$. Thus, by (\ref{philbn}),
$\pl(\p(a))\in f(\s)$. Using ($\LAM$2), we obtain that $a\in
f(\s)$. Therefore $f(\s)\in\l_A^g(a)=F$. So, $\s\in f\inv(F)$. We
have shown that $G\sbe f\inv(F)$. Then
$\int(G)\sbe\cl(\int(f\inv(F)))$ and thus $G\sbe\ps(F)$.
Conversely, let $\s\in\cl(f\inv(\int(F)))$. Suppose that
$\p(a)\nin\s$. Then $\s\in\int(\l_B^g(\p(a^*)))$ (see
(\ref{inthal})).
 Hence $\int(\l_B^g(\p(a^*)))\cap
f\inv(\int(F))\nes$. Thus there exists
$\s\ap\in\int(\l_B^g(\p(a^*)))$ such that $f(\s\ap)\in\int(F)$.
Then, by (\ref{inthal}), $a^*\nin f(\s\ap)$. Since
$\p(a^*)\in\s\ap$, (\ref{philbn}) and ($\LAM$2) imply that $a^*\in
f(\s\ap)$, a contradiction. Hence $\p(a)\in\s$, i.e. $\s\in G$.
So, $\ps(F)\sbe G$. Therefore, $\ps(F)=G$. This shows that
\begin{equation}\label{psiphininnk}
\l_B^g\circ\p=\ps\circ \l_A^g.
\end{equation}
Hence, $\l$ is a natural isomorphism.

We have proved that $\SKLC$ and $\SKAL$ are dually equivalent
categories.
\sqs

\begin{defi}\label{categn}{\rm (Fedorchuk \cite{F})}
\rm We will denote by $\SHC$ the category of all compact Hausdorff
spaces and all quasi-open maps between them.

Let $\SAC$ be the category whose objects are all complete normal
contact algebras and whose morphisms $\p:(A,C)\lra (B,C\ap)$ are
all complete Boolean homomorphisms $\p:A\lra B$ satisfying the
following condition:

\smallskip

\noindent (F1) For all $a,b\in A$,  $\p(a)C\ap\p(b)$ implies
$aCb$.

\medskip

It is easy to see that in this  way we have defined categories.
\end{defi}

\begin{theorem}\label{dcompn}{\rm (Fedorchuk \cite{F})}
The categories $\SHC$ and $\SAC$ are dually equi\-valent.
\end{theorem}

\doc It follows from Theorem \ref{maintheoremnk}, \ref{conaln} and \ref{semsk}(b).
\sqs

We will now obtain  one more generalization of Theorem
\ref{dcompn}.

\begin{defi}\label{lcatn}
\rm Let $\SLC$ be the category of all locally compact Hausdorff
spaces and all skeletal perfect maps between them. Note that, by
\ref{semsk}(b), the morphisms of the category $\SLC$ are precisely
the  quasi-open  perfect maps.

Let $\SAL$ be the category whose objects are all complete local
contact algebras (see \ref{locono}) and whose morphisms are all
$\SKAL$-morphisms $\p:(A,\rho,\BBBB)\lra (B,\eta,\BBBB\ap)$
satisfying the following condition:

\smallskip

\noindent(L3) $a\in\BBBB$ implies $\p(a)\in\BBBB\ap$.

\smallskip

It is easy to see that in this  way we have defined categories.
Obviously, $\SLC$ (resp., $\SAL$) is a (non-full) subcategory of
the category $\SKLC$ (resp., $\SKAL$).
\end{defi}


\begin{theorem}\label{maintheoremn}
The categories $\SLC$ and $\SAL$ are dually equivalent.
\end{theorem}

\doc We will show that the restrictions
$\Psi^a_p:\SAL\lra\SLC$ and $\Psi^t_p:\SLC\lra\SAL$ of the
contravariant functors $\Psi^a$ and $\Psi^t$ defined in the proof
of Theorem \ref{maintheoremnk} are the desired duality functors.

Let $f\in\SLC((X,\tau),(Y,\tau\ap))$. Since $f$ is a perfect map,
we obtain that $\p=\Psi^t_p(f)$ satisfies  condition (L3) (using
\cite[Theorem 3.7.2]{E}). Hence, $\p$ is well defined. Therefore
the contravariant functor $\Psi^t_p:\SLC\lra\SAL$ is well defined.

Let  $\p\in\SAL((A,\rho,\BBBB),(B,\eta,\BBBB\ap))$ and set
$f=\Psi^a_p(\p)$, i.e. $f:\Psi^a_p(B,\eta,\BBBB\ap)\lra
\Psi^a_p(A,\rho,\BBBB)$. Put $C=C_\rho$ and $C\ap=C_\eta$ (see
\ref{Alexprn} for the notations). Then, by \ref{Alexprn}, $(A,C)$
and $(B,C\ap)$ are CNCA's.
 Denote by $\p_c$ the map $\p$ regarded as a function of $(A,C)$ to
 $(B,C\ap)$.
 We will show that $\p_c$ satisfies
condition (F1)  in \ref{categn}.

 For verifying  condition (F1), let
$a,b\in A$ and let $\p_c(a) C\ap \p_c(b)$. Then either
$\p_c(a)\eta\p_c(b)$ or $\p_c(a),\p_c(b)\nin\BBBB\ap$. If
$\p_c(a)\eta\p_c(b)$ then, by (L1), $a\rho b$; hence $aCb$. If
$\p_c(a),\p_c(b)\nin\BBBB\ap$ then, by (L3), $a,b\nin\BBBB$. Hence
$a C b$. So, (F1) is verified.  Therefore,
\begin{equation}\label{M1i3n}
\p_c:(A,C_\rho)\lra (B,C_\eta) \mbox{ satisfies condition (F1).}
\end{equation}

Set $X=\Psi^a(A,C,A)$ and $Y=\Psi^a(B,C\ap,B)$ (see
(\ref{phiapcn})). Then $X$ and $Y$ are compact Hausdorff spaces.
Let $f_c=\Psi^a(\p_c)$, i.e.
\begin{equation}\label{fden}
f_c: Y\lra X \mbox{ is defined by } f_c(\s_u)=\s_{\p_c\inv(u)},
\mbox{ for every }u\in \Ult(B).
\end{equation}
 Then, by (\ref{fcontk}), (\ref{skm}) and \ref{semsk}(b), we obtain that
\begin{equation}\label{fcontn}
f_c:Y\lra X \mbox{ is a quasi-open map.}
\end{equation}

We will regard three cases now.

\medskip

\noindent(a)~ Let $1_A\nin\BBBB$ and $1_B\nin\BBBB\ap$. Then
$\Psi^a_p(B,\eta,\BBBB\ap)=L_B=Y\stm\{\s_\infty^B\}$ and
$\Psi^a_p(A,\rho,\BBBB))=L_A=X\stm\{\s_\infty^A\}$ (see
\ref{neogrn} and (\ref{L})).

 We will show
that $f_c\inv(\s_\infty^A)=\{\s_\infty^B\}$ (see \ref{neogrn} for
the notations). We first prove that
$f_c(\s_\infty^B)=\s_\infty^A$. Let $u\in \Ult(B)$ be such that
$u\subset\s_\infty^B$ and $\s_\infty^B=\s_u$ (see
\ref{conclustth}). Then $f_c(\s_\infty^B)=\s_{\p_c\inv(u)}$. We
will show that $\p_c\inv(u)\subset \s_\infty^A$. Indeed, let
$a\in\p_c\inv(u)$. Then $\p_c(a)\in u\subset B\stm\BBBB\ap$. Hence
$\p_c(a)\nin\BBBB\ap$. Thus, by (L3), $a\nin \BBBB$. So,
$\p_c\inv(u)\subset A\stm\BBBB=\s_\infty^A$ (see \ref{neogrn}).
Then, by \ref{neogrn} and \ref{uniqult},
$\s_\infty^A=\s_{\p_c\inv(u)}$. Therefore,
$f_c(\s_\infty^B)=\s_\infty^A$. Since $L_A$ and $L_B$ consist of
bounded clusters (see (\ref{L})), (\ref{sbounn}) implies that
$f_c(L_B)\sbe L_A$. Therefore,
$f_c\inv(\s_\infty^A)=\{\s_\infty^B\}$. This shows that
$f_c\inv(L_A)=L_B$. Since $f_c$ is a  perfect map, we obtain (by
\cite[Proposition 3.7.4]{E}) that
\begin{equation}\label{fldefn}
(f_c)_{L_A}:L_B\lra L_A \mbox{ is a perfect map.}
\end{equation}
Obviously, $f$ is the restriction of $f_c$ to $L_B$. Hence
$f=(f_c)_{L_A}$, i.e. $f$ is a perfect map. Since $f$ is a
skeletal map (by (\ref{skm})), \ref{semsk}(b) implies that
\begin{equation}\label{fln}
f \mbox{ is a quasi-open perfect map.}
\end{equation}

\medskip

\noindent(b)~ Let $1_A\nin\BBBB$ and $1_B\in\BBBB\ap$. Then
$C\ap=\eta$, $\Psi^a_p(A,\rho,\BBBB)=X\stm\{\s_\infty^A\}=L_A$ and
$\Psi^a_p(B,\eta,\BBBB\ap)=Y$. Thus (\ref{sbounn}) implies that
$f_c(Y)\subset L_A$. Therefore, the restriction $f:Y\lra L_A$ of
$f_c$ is a perfect map. Since $f$ is skeletal (by (\ref{skm})), we
obtain, using \ref{semsk}(b), that $f$ is quasi-open. Therefore,
\begin{equation}\label{fapn}
f
:\Psi^a_p(B,\eta,\BBBB\ap)\lra\Psi^a_p(A,\rho,\BBBB)
\end{equation}
is a quasi-open perfect map.

\medskip

\noindent(c)~ Let $1_A\in\BBBB$. Then, by (L3), $1_B\in\BBBB\ap$.
Hence $C=\rho$, $C\ap=\eta$, $\Psi^a_p(B,\eta,\BBBB\ap)=Y$,
$\Psi^a_p(A,\rho,\BBBB)=X$. Thus $f=f_c$. Hence, by
(\ref{fcontk}), (\ref{skm}) and \ref{semsk}(b), $f:Y\lra X$ is a
quasi-open perfect map.

\medskip

We have regarded all possible cases.
 Therefore, $\Psi^a_p$ is well defined on the objects
and morphisms of the category $\SAL$.

Note that, using (\ref{hapisomn}), we obtain
that  $\l^g_B$ is a $\SAL$-isomorphism. The rest follows from
Theorem \ref{maintheoremnk}.
\sqs

\begin{defi}\label{lcat}
\rm Let $\OLC$ be the category of all locally compact Hausdorff
spaces and all open maps between them.

Let $\OAL$ be the category whose objects are all complete local
contact algebras
and whose morphisms are all $\SKAL$-morphisms
$\p:(A,\rho,\BBBB)\lra (B,\eta,\BBBB\ap)$ satisfying the following
condition:

\smallskip

\noindent (LO) $\fa a\in A$ and $\fa b\in \BBBB\ap$,
$\p_\LAM(b)\rho a$ implies
 $b\eta\p(a)$.

\smallskip

It is easy to see that in this  way we have defined categories.
Obviously, $\OAL$ (resp., $\OLC$) is a (non-full) subcategory of
the category $\SKAL$ (resp., $\SKLC$).
\end{defi}


\begin{theorem}\label{maintheorem}
The categories $\OLC$ and $\OAL$ are dually equivalent.
\end{theorem}

\doc We will show that the restrictions
$\Psi^a_o:\OAL\lra\OLC$ and $\Psi^t_o:\OLC\lra\OAL$ of the
contravariant functors $\Psi^a$ and $\Psi^t$ defined in the proof
of Theorem \ref{maintheoremnk} are the desired duality functors.

Let $f\in\OLC((X,\tau),(Y,\tau\ap))$. Set $\p= \Psi^t_o(f)$. Then,
since $f$ is an open map, \cite[1.4.C]{E}  implies that for every
$F\in RC(Y)$,
$f\inv(F)=f\inv(\cl(\int(F)))=\cl(f\inv(\int(F)))=\p(F)$ (see
(\ref{defpsiapnk})). Hence,
\begin{equation}\label{defpsiap}
\Psi^t_o(f):\Psi^t_o(Y,\tau\ap)\lra\Psi^t_o(X,\tau) \mbox{ is
defined by } \Psi^t(f)(F)=f\inv(F),
\end{equation}
for all $F\in\Psi^t_o(Y,\tau\ap)$.  Further, by the proof of
Theorem \ref{maintheoremnk}, $\p$ is an $\SKAL$-morphism. We will
show that $\p$ satisfies condition (LO). We have that
$\pl:RC(X)\lra RC(Y)$ is defined, according to (\ref{leftadjnk})
and (\ref{psilnk}), by the formula $\pl(F)=\cl(f(F))$, for every
$F\in RC(X)$. So, let $F\in RC(Y)$, $G\in CR(X)$ and
$F\rho_Y\p_\LAM (G)$; then $F\cap f(G)\nes$ and hence
$f\inv(F)\cap G\nes$; therefore, $\p(F)\rho_X G$. So, the axiom
(LO) is fulfilled.
 Hence, $\Psi^t_o(f)$ is an $\OAL$-morphism. Therefore, the
contravariant  functor  $\Psi^t_o$ is well defined.

Let  $\p\in\OAL((A,\rho,\BBBB),(B,\eta,\BBBB\ap))$. Put $C=C_\rho$
and $C\ap=C_\eta$ (see \ref{Alexprn} for the notations). Then, by
\ref{Alexprn}, $(A,C)$ and $(B,C\ap)$ are CNCA's.

Set $X=\Psi^a_o(A,\rho,\BBBB)$, $Y=\Psi^a_o(B,\eta,\BBBB\ap)$  and
$f=\Psi^a_o(\p)$. Then, by the proof of Theorem
\ref{maintheoremnk}, $f:Y\lra X$ is a continuous skeletal map. We
are now going to show that $f$ is an open map. By (\ref{eel}), it
is enough to prove that, for every $b\in\BBBB\ap$,
$f(\int_{Y}(\l_B(b)))$ is an open subset of $X$ (note that
$\l_B(b)=\l_B^g(b)$ because $b\in\BBBB\ap$).

So, let $b\in \BBBB\ap$. Let $\s\in f(\int_{Y}(\l_B(b)))$. Then
there exists $\s\ap\in\int_{Y}(\l_B(b))$ such that $\s=f(\s\ap)$.
By (\ref{inthal}), $b^*\not\in\s\ap$. Then \ref{cluendcor} implies
that there exists $c_1\in B$ such that $b^*\ll_{C\ap} c_1^*$ and
$c_1^*\not\in \s\ap$. Since $\s\ap$ is a bounded cluster in
$(B,C\ap)$, (\ref{bstar}) implies that there exists
$c_2\in\BBBB\ap$ such that $c_2^*\nin\s\ap$. Put $b_1=c_1\we c_2$.
Then $b_1\in\BBBB\ap\cap\s\ap$ (by (\ref{b1b2})), $b_1^*\nin\s\ap$
(by (\ref{b1b2})) and $b^*\ll_{C\ap} b_1^*$ (by ($\ll$3) (see
\ref{conalg})). Thus $b_1\ll_{C\ap} b$. Therefore, by
(\ref{inthal}) and (\ref{abipn}), $\s\ap\in\int_{Y}(\l_B(b_1))\sbe
\l_B(b_1)\sbe\int_{Y}(\l_B(b))$. By \ref{conclustth}, there exists
$u\in \Ult(B)$ such that $b_1\in u\sbe\s\ap$ and $\s\ap=\s_u$. Put
$a=\p_\LAM (b_1)$. Then, by (\ref{philbn}), $a\in f(\s\ap)=\s$.
Suppose that $a^*\in\s$. We will show that this implies that
$\p(a^*)\in\s\ap$. Indeed, suppose that $\p(a^*)\nin\s\ap$. Then
there exists $c_3\in u$ such that $\p(a^*)(-C\ap)c_3$. Set
$b_2=c_2\we c_3$. Then $b_2\in u\cap\BBBB\ap$ and
$\p(a^*)(-C\ap)b_2$.  Since $C\ap=C_\eta$, we obtain, by
\ref{Alexprn}, that $\p(a^*)(-\eta)b_2$. Using condition (LO), we
get that $a^*(-\rho)\p_\LAM(b_2)$. Since $\p_\LAM(b_2)\in\BBBB$
(by (L2)), we obtain that $a^*(-C)\p_\LAM(b_2)$ (see again
\ref{Alexprn}).
 By ($\LAM$1), $\p(\p_\LAM(b_2))\ge b_2$;  thus
$\p(\p_\LAM(b_2))\in u$. Hence $\p_\LAM(b_2)\in\p\inv(u)$. Since
$\s=f(\s\ap)=\s_{\p\inv(u)}$ and $a^*\in\s$, we have that $a^*Cc$,
for every $c\in\p\inv(u)$. Therefore $a^*C\p_\LAM(b_2)$,  a
contradiction. Hence, $\p(a^*)\in\s\ap$, i.e. $(\p(\p_\LAM
(b_1)))^*\in\s\ap$. Since, by ($\LAM$1),  $b_1^*\ge (\p(\p_\LAM
(b_1)))^*$, we obtain that $b_1^*\in\s\ap$,  a contradiction.
Thus, $a^*\not\in\s$. Then, using (\ref{inthal}), (\ref{abipn})
and (\ref{raven2n}), we obtain that $\s\in\int_{X}(\l_A(a))\sbe
\l_A(a)=\l_A(\p_\LAM (b_1))=f(\l_B(b_1))\sbe
f(\int_{Y}(\l_B(b)))$. Therefore, $f(\int_{Y}(\l_B(b)))$ is an
open set in $X$. Thus, $f$ is an open map. Hence $\Psi^a_o$ is
well defined.

Further, note that, using (\ref{hapisomn}), it is easy to see that
$\l^g_B$ is an $\OAL$-isomorphism. The rest follows from Theorem
\ref{maintheoremnk}.
\sqs

\begin{defi}\label{categ}
\rm We will denote by $\OHC$ the category of all compact Hausdorff
spaces and all open maps between them.

Let $\OAC$ be the category whose objects are all complete normal
contact algebras and whose morphisms are all $\SAC$-morphisms
$\p:(A,C)\lra (B,C\ap)$ satisfying the following condition:

\smallskip

\noindent (CO) For all $a\in A$ and all $b\in B$, $aC\p_\LAM (b)$
implies $\p(a)C\ap b$ (see \ref{ladj}  for $\p_\LAM $).

\medskip

It is easy to see that in this  way we have defined categories.
The category $\OAC$ (resp., $\OHC$) is a (non-full) subcategory of
the category $\SAC$ (resp., $\SHC$).
\end{defi}

\begin{theorem}\label{dcomp}
The categories $\OHC$ and $\OAC$ are dually equivalent.
\end{theorem}

\doc It follows from Theorem \ref{maintheorem} and \ref{conaln}.
\sqs

\begin{defi}\label{lcatp}
\rm Let $\OPLC$ be the category of all locally compact Hausdorff
spaces and all open  perfect maps between them.

Let $\OPAL$ be the category whose objects are all complete local
contact algebras (see \ref{locono}) and whose morphisms are all
$\SAL$-morphisms  satisfying  condition (LO).


It is easy to see that in this  way we have defined categories.
Obviously, $\OPAL$ (resp., $\OPLC$) is a subcategory of the
category $\SAL$ (resp., $\SLC$).
\end{defi}

\begin{theorem}\label{maintheoremp}
The categories $\OPLC$ and $\OPAL$ are dually equivalent.
\end{theorem}

\doc It follows from Theorems \ref{maintheoremn} and
\ref{maintheorem}.
\sqs

Note that since the morphisms of the category $\OPLC$ are closed
maps, in the definition of the category $\OPAL$ (see \ref{lcatp})
we can substitute condition (LO)   for the following one:

\smallskip

\noindent (LO') $\fa a\in A$ and $\fa b\in B$, $a\rho\p_\LAM(b)$
implies $\p(a)\eta b$.

\section{Connected Spaces}


\begin{notas}\label{convcatn}
\rm If $\K$ is a category whose objects form a subclass of the
class of all topological spaces (resp., contact algebras) then we
will denote by $\KCon$ the full subcategory of $\K$ whose objects
are all $``$connected" $\K$-objects, where $``$connected" is
understood in the usual sense when the objects of $\K$ are
topological spaces and in the sense of \ref{conalg} (see the
condition (CON) there) when the objects of $\K$ are contact
algebras. For example, we  denote by:


\noindent $\bullet$ $\SLCC$ the full subcategory of the category
$\SLC$ having as objects  all connected locally compact Hausdorff
spaces;


\noindent $\bullet$ $\SALC$ the full subcategory of the category
$\SAL$ having as objects all connected CLCA's.

\end{notas}

\begin{theorem}\label{conthen}
The categories $\SLCC$ and $\SALC$ are dually equivalent; in
particular, the categories $\SHCC$ and $\SACC$ are dually
equivalent.
%
\end{theorem}

\doc It follows immediately  from \ref{confact},
Theorem \ref{maintheoremn} and Theorem \ref{dcompn}. \sqs

%
%
%
%

\begin{theorem}\label{conthe}
The categories $\OPLCC$ and $\OPALC$ are dually equ\-i\-valent; in
particular, the categories $\OHCC$ and $\OACC$ are dually
equivalent.
\end{theorem}

\doc It follows immediately from \ref{confact},
Theorem \ref{dcomp} and Theorem \ref{maintheoremp}. \sqs

Analogously one can formulate and prove the connected versions of
the theorems Theorem \ref{maintheoremnk} and Theorem
\ref{maintheorem}.


\section{Equivalence Theorems}

\begin{defi}\label{feddef}{\rm (\cite{F})}
\rm
 Let $\ESHC$ be the category whose objects are all complete
normal contact algebras and whose morphisms
$\psi:(A,C)\lra(B,C\ap)$  are all functions $\psi:A\lra B$
satisfying the following conditions:

\medskip

\noindent(EF1) for every $a\in A$, $\psi(a)=0$ iff $a=0$;\\
(EF2) $\psi$ preserves all joins;\\
(EF3) if $a\in A$, $b\in B$ and $b\le\psi(a)$ then there exists
$c\in A$ such that $c\le a$ and $\psi(c)=b$;\\
(EF4) for every $a,b\in A$, $aCb$ implies that
$\psi(a)C\ap\psi(b)$.
\end{defi}

In \cite{F}, V. V. Fedorchuk proved the following theorem:

\begin{theorem}\label{fedeth}{\rm (\cite{F})}
The categories\/ $\SHC$ and\/ $\ESHC$ are equivalent.
\end{theorem}

We will now present a generalization of this theorem.

\begin{defi}\label{lcedef}
\rm
 Let $\ESKLC$ be the category whose objects are all complete
local contact algebras and whose morphisms
$\psi:(A,\rho,\BBBB)\lra (B,\eta,\BBBB\ap)$ are all functions
$\psi:A\lra B$  satisfying conditions (EF1)-(EF3) (see Definition
\ref{feddef}) and the following two constraints:

\medskip

\noindent(EL4) for every $a,b\in A$, $a\rho b$ implies that
$\psi(a)\eta\psi(b)$;\\
(EL5) if $a\in\BBBB$ then $\psi(a)\in\BBBB\ap$.
\end{defi}

The proof of the following theorem is similar to that of Theorem
\ref{fedeth}.

\begin{theorem}\label{lceth}
The categories\/ $\SKLC$ and\/ $\ESKLC$ are equivalent.
\end{theorem}

\doc Since the categories $\SKLC$ and $\SKAL$ are dually
equivalent (by Theorem \ref{maintheoremnk}), it is enough to show
that the categories $\ESKLC$ and $\SKAL$ are dually equivalent.

Let us define a contravariant functor $D_p: \ESKLC\lra\SKAL$. Let
$D_p$  be the identity on the objects of the category $\ESKLC$ and
let, for every $\psi\in\ESKLC((A,\rho,\BBBB),(B,\eta,\BBBB\ap))$,
$D_p(\psi)=\psi_P$, where $\psi_P$ is the right adjoint of $\psi$
(see \ref{ladj} and (EF2)). Setting $\p=\psi_P$, we have to show
that $$\p\in\SKAL((B,\eta,\BBBB\ap),(A,\rho,\BBBB)).$$

As it is proved in \cite{F}, $\p$ is a complete Boolean
homomorphism. For completeness of our exposition, we will present
here the Fedorchuk's proof. Note first that $\psi=\pl$. By
\ref{ladj}, $\p$ preserves all meets in $B$. Since, by (EF1),
$\psi(0)=0$, we have that $\p(0)=\p(\psi(0))$; if $\p(0)>0$ then,
by (EF1) and \ref{ladj}, $0=\psi(0)=\psi(\p(\psi(0)))>0$, a
contradiction. Hence $\p(0)=0$. Further, since $\psi(1)\le 1\iff
1\le\p(1)$, we get that $\p(1)=1$. Finally, $\p(b^*)=(\p(b))^*$,
for every $b\in B$. Indeed, let $b\in B$. Set $a=\p(b)\we\p(b^*)$.
Then, by \ref{ladj}, $\psi(a)\le\psi(\p(b))\we\psi(\p(b^*))\le
b\we b^*=0$. Hence $\psi(a)=0$. Therefore, by (EF1), $a=0$, i.e.
$\p(b)\we\p(b^*)=0$. Set now $c=\p(b)\vee\p(b^*)$ and suppose that
$c<1$. Then $c^*\neq 0$. Since $0=c^*\we
c=(c^*\we\p(b))\vee(c^*\we\p(b^*))$, we have that
$c^*\we\p(b)=0=c^*\we\p(b^*)$. By (EF1), $\psi(c^*)\neq 0$.
Obviously, $\psi(c^*)=(\psi(c^*)\we b)\vee(\psi(c^*)\we b^*)$.
Therefore, at least one of the elements $\psi(c^*)\we b$ and
$\psi(c^*)\we b^*$ is different from $0$. Let $\psi(c^*)\we b\neq
0$.  By (EF3), the inequality $\psi(c^*)\we b\le\psi(c^*)$ implies
that there exists $d\in A$ such that $d\le c^*$ and
$\psi(d)=\psi(c^*)\we b$. Since $\psi(d)\neq 0$, we get, by (EF1),
that $d\neq 0$. Further, $\psi(d)\le b$ implies that $d\le\p(b)$.
Then $d\le c^*\we\p(b)=0$, i.e. $d=0$, a contradiction.
Analogously, we obtain a contradiction if $\psi(c^*)\we b^*\neq
0$. So, $c=1$, i.e. $\p(b)\vee\p(b^*)=1$. Hence, we have proved
that $\p(b^*)=(\p(b))^*$. All this shows that $\p$ is a complete
Boolean homomorphism.

Since conditions (L1) and (EL1) in \ref{lcatnk} are equivalent and
$\psi=\pl$, (EL4) implies that $\p$ satisfies condition (L1).
Obviously, (EL5) implies that $\p$ satisfies condition (L2) in
\ref{lcatnk}. So, $\p$ is a $\SKAL$-morphism. Now, from
$D_p(id)=id$ and the formula
$(\psi_2\circ\psi_1)_P=(\psi_1)_P\circ(\psi_2)_P$, we obtain that
$D_p$ is a contravariant functor.

Let us define a contravariant functor $D_l: \SKAL\lra\ESKLC$. Let
$D_l$  be the identity on the objects of the category $\SKAL$ and
let, for every $\p\in\SKAL((A,\rho,\BBBB),(B,\eta,\BBBB\ap))$,
$D_l(\p)=\pl$, where $\pl$ is the left adjoint of $\p$ (see
\ref{ladj}). Setting $\psi=\pl$, we have to show that
$\psi\in\ESKLC((B,\eta,\BBBB\ap),(A,\rho,\BBBB)).$

Since $0\le\p(0)$ implies that $\psi(0)\le 0$, we get that
$\psi(0)=0$. If $\psi(b)=0$ then $\psi(b)\le 0$ and hence
$b\le\p(0)=0$, i.e. $b=0$. Therefore, $\psi$ satisfies condition
(EF1). Further, conditions (EF2), (EL4) and (EL5) are clearly
satisfied by $\psi$. Finally, let $a\le\psi(b)$. Set
$c=b\we\p(a)$. Then $c\le b$ and, by \ref{L2rave}(b),
$\psi(c)=a\we\psi(b)=a$. Therefore, $\psi$ satisfies condition
(EF3). So, $\psi$ is an $\ESKLC$-morphism. Now, it is clear that
$D_l$ is a contravariant functor. Since the compositions of $D_p$
and $D_l$ are the identity functors, we get that $D_p$ is a
duality. Put now $\Phi^a=\Psi^a\circ D_p$ and
$\Phi^t=D_l\circ\Psi^t$. Then $\Phi^a:\ESKLC\lra\SKLC$ and
$\Phi^t:\SKLC\lra\ESKLC$ are the required equivalences.
 \sqs

\begin{defi}\label{lcskpedef}
\rm
 Let $\ESLC$ be the category whose objects are all complete
local contact algebras (see \ref{locono}) and whose morphisms are
all $\ESKLC$-morphisms $\psi:(A,\rho,\BBBB)\lra (B,\eta,\BBBB\ap)$
satisfying the following condition:

\medskip

\noindent(EL6) if $b\in\BBBB\ap$ then $\psi_P(b)\in\BBBB$ (where
$\psi_P$ is the right adjoint of $\psi$ (see \ref{ladj})).
\end{defi}

\begin{theorem}\label{lcskpeth}
The categories\/ $\SLC$ and\/ $\ESLC$ are equivalent.
\end{theorem}

\doc Using Theorem \ref{maintheoremn}, it is enough to show that the categories
$\SAL$ and $\ESLC$ are dually equivalent. We will show that the
restriction of the contravariant functor $D_p$ (defined in the
proof of Theorem \ref{lceth}) to the category $\ESLC$ is the
required duality functor.

Let $\psi\in\ESLC((A,\rho,\BBBB),(B,\eta,\BBBB\ap))$. Then, by
(EL6), $\psi_P$ satisfies condition (L3) in \ref{lcatn}. Hence, by
the proof of Theorem \ref{lceth}, $D_p(\psi)$ is a
$\SAL$-morphism. Further, let us regard the restriction of the
contravariant functor $D_l$ (defined in the proof of Theorem
\ref{lceth}) to the category $\SAL$. If $\p$ is a $\SAL$-morphism
then, by (L3), $\pl$ satisfies condition (EL6). Hence $D_l(\p)$ is
an $\ESLC$-morphism. Therefore, $D_p$ is a duality. \sqs

\begin{defi}\label{lcoedef}
\rm
 Let $\EOLC$ be the category whose objects are all complete
local contact algebras and whose morphisms are all
$\ESKLC$-morphisms $\psi:(A,\rho,\BBBB)\lra (B,\eta,\BBBB\ap)$
satisfying the following condition:

\medskip

\noindent(EL7) if $b\in B$, $a\in\BBBB$ and $\psi(a)\eta b$ then
$a\rho\psi_P(b)$ (where $\psi_P$ is the right adjoint of $\psi$
(see \ref{ladj})).
\end{defi}

\begin{theorem}\label{lcoeth}
The categories\/ $\OLC$ and\/ $\EOLC$ are equivalent.
\end{theorem}

\doc It is clear that if $\psi$ satisfies condition (EL7) then
$\psi_P$ satisfies condition (LO) in \ref{lcat} and if $\p$
satisfies condition (LO) then $\pl$ satisfies (EL7). Now, using
Theorem \ref{maintheorem}, we argue as in the proof of Theorem
\ref{lcskpeth}. \sqs

\begin{defi}\label{coedef}
\rm
 Let $\EOPC$ be the category whose objects are all complete
normal contact algebras and whose morphisms are all
$\ESHC$-morphisms $\psi:(A,C)\lra(B,C\ap)$ satisfying the
following condition:

\medskip

\noindent(EC7) if $a\in A$, $b\in B$ and $\psi(a)\, C\ap\, b $
then $a\, C\, \psi_P(b)$ (where $\psi_P$ is the right adjoint of
$\psi$ (see \ref{ladj})).
\end{defi}

\begin{theorem}\label{coeth}
The categories\/ $\OHC$ and\/ $\EOPC$ are equivalent.
\end{theorem}

\doc It follows directly from Theorem \ref{lcoeth}. \sqs

\begin{defi}\label{lcopedef}
\rm
 Let $\EOPLC$ be the category whose objects are  all complete
local contact algebras and whose morphisms are all
$\ESLC$-morphisms satisfying condition (EL7).
\end{defi}

\begin{theorem}\label{lcopeth}
The categories\/ $\OPLC$ and\/ $\EOPLC$ are equivalent.
\end{theorem}

\doc It follows from the proofs of \ref{lcskpeth} and
\ref{lcoeth}. \sqs

\begin{rem}\label{Ro}
\rm A great part  of our Theorem \ref{lceth} is formulated (in
another form) and proved
 in Roeper's paper \cite{R}. Let us
state precisely what is done there (using our notations). Roeper
defines the notion of {\em mereological mapping}\/: such is any
function $\psi:B\lra A$, where $A$ and $B$ are complete Boolean
algebras, which satisfies the following conditions: (i)
$\psi(b)=0$ iff $b=0$; (ii) $a\le b$ implies $\psi(a)\le\psi(b)$;
(iii) if $0\neq a\le\psi(b)$, where $b\in B$ and $a\in A$, then
there exists $b\ap\in B$ such that $0\neq b\ap\le b$ and
$\psi(b\ap)\le a$. It is shown that any mereological mapping
preserves all joins in $B$. Further, a mapping $\psi$ of a CLCA
$(B,\eta,\BBBB\ap)$ to another CLCA $(A,\rho,\BBBB)$ is called:
(a) {\em continuous}\/ if $a\eta b$ implies $\psi(a)\rho\psi(b)$,
and (b) {\em bounded}\/ if $\psi(b)\in\BBBB$ when $b\in\BBBB\ap$.
It is shown that every continuous and bounded mereological mapping
$\psi:(B,\eta,\BBBB\ap)\lra(A,\rho,\BBBB)$ generates a function
$f_\psi:\Psi^a(B,\eta,\BBBB\ap)\lra\Psi^a(A,\rho,\BBBB)$, defined
by the formula $f_\psi(\s_u)=\s_{\psi(u)}$, for every
$u\in\Ult(B)$; the function $f_\psi$ is continuous (in topological
sense) and is such that $\cl(f_\psi(F))$ is regular closed when
$F$ is regular closed. It is proved that if
$f:\Psi^a(B,\eta,\BBBB\ap)\lra\Psi^a(A,\rho,\BBBB)$ is a
continuous function such that $\cl(f(F))$ is regular closed when
$F$ is regular closed then there exists a continuous and bounded
mereological function $\psi:(B,\eta,\BBBB\ap)\lra (A,\rho,\BBBB)$
such that $f=f_\psi$. Finally, a mereological function
$\psi:(B,\eta,\BBBB\ap)\lra (A,\rho,\BBBB)$ is called {\em
topological}\/ if $\psi(1_B)=1_A$, $\psi(a)\rho\psi(b)$ iff $a\eta
b$, and $\psi(b)\in\BBBB$ iff $b\in\BBBB\ap$; it is shown that if
$\psi$ is topological then $f_\psi$ is a homeomorphism and if
$f:\Psi^a(B,\eta,\BBBB\ap)\lra\Psi^a(A,\rho,\BBBB)$ is a
homeomorphism then there exists a topological function
$\psi:(B,\eta,\BBBB\ap)\lra (A,\rho,\BBBB)$ such that $f=f_\psi$.

It is easy to see that a function $\psi:B\lra A$ is mereological
iff it satisfies conditions (EF1)-(EF3) (see Definition
\ref{feddef}); $\psi$ is continuous (respectively, bounded) iff it
satisfies condition  (EL4) (respectively, (EL5)).
 Further, Lemma
\ref{mrro} shows that a continuous map $f:X\lra Y$ satisfies
Roeper's condition $``\cl(f(F))\in RC(Y)$ when $F\in RC(X)$" iff
$f$ is a skeletal map. Therefore, our covariant functor
$\Phi^a:\ESKLC\lra\SKLC$ (see the proof of Theorem \ref{lceth})
was defined in \cite{R} in another but equivalent form and  it was
shown there that $\Phi^a$ is full and isomorphism-dense; however,
in \cite{R} it was not shown that $\Phi^a$ is faithful.
\end{rem}


\baselineskip = 0.75\normalbaselineskip

\end{document}